\documentstyle[amscd,amssymb]{amsart}
\newcommand{\nc}{\newcommand}

\nc{\ad}{{\mbox{\bf{ad}}}}
\nc{\AJ}{{\operatorname{aj}}}
\nc{\Aut}{{\operatorname{Aut}}}
\nc{\Bls}{{{\cal B}ls}}
\nc{\Boxtimes}{{\fbox{$\times$}}}
\nc{\blt}{{\bullet}}
\nc{\bSt}{{\mbox{\bf{St}}}}
\nc{\card}{{\operatorname{card}}}
\nc{\Cch}{{\check{C}}}
\nc{\cd}{{\operatorname{cd}}}
\nc{\Ch}{{\operatorname{Ch}}}
\nc{\chara}{{\operatorname{char}}}
\nc{\CHom}{{\cal{H}om}}
\nc{\Coker}{{\operatorname{Coker}}}
\nc{\codim}{{\operatorname{codim}}}
\nc{\Cone}{{\operatorname{Cone}}}
\nc{\cSgn}{{\cal{S}gn}}
\nc{\depth}{{\operatorname{depth}}}
\nc{\dirlim}{{\underset{\rightarrow}{\operatorname{lim}}}}
\nc{\dotbox}{{\overset{\bullet}{\boxtimes}}}
\nc{\dotimes}{{\overset{\bullet}{\otimes}}}
\nc{\Ed}{{\operatorname{Edge}}}
\nc{\emp}{{\emptyset}}
\nc{\Ext}{{\operatorname{Ext}}}
\nc{\Fac}{{\cal{F}ac}}
\nc{\Fun}{{\operatorname{F}}}
\nc{\FS}{{\cal{FS}}}
\nc{\Hom}{{\operatorname{Hom}}}
\nc{\had}{{{\hat{\mbox{\bf{ad}}}}}}
\nc{\hgt}{{\operatorname{ht}}}
\nc{\Id}{{\operatorname{Id}}}
\nc{\id}{{\operatorname{id}}}
\nc{\Ima}{{\operatorname{Im}}}
\nc{\ind}{{\operatorname{ind}}}
\nc{\Ind}{{\operatorname{Ind}}}
\nc{\infi}{{\operatorname{inf}}}
\nc{\infh}{{\frac{\infty}{2}}}
\nc{\invlim}{{\underset{\leftarrow}{\operatorname{lim}}}}
\nc{\Jac}{{{\cal J}ac}}
\nc{\Ker}{{\operatorname{Ker}}}
\nc{\lcm}{{\operatorname{lcm}}}
\nc{\length}{{\operatorname{length}}}
\nc{\Locsys}{{{\cal L}ocsys}}
\nc{\Map}{{{\cal M}ap}}
\nc{\modul}{{\operatorname{mod}}}
\nc{\Mor}{{\operatorname{Mor}}}
\nc{\MS}{{\cal{MS}}}
\nc{\Ob}{{\operatorname{Ob}}}
\nc{\opp}{{\operatorname{opp}}}
\nc{\Or}{{{\cal O}r}}
\nc{\Ord}{{{\cal O}rd}}
\nc{\Part}{{{\cal P}art}}
\nc{\PGL}{{\operatorname{PGL}}}
\nc{\Pic}{{\operatorname{Pic}}}
\nc{\Rep}{{{\cal{R}}ep}}
\nc{\rk}{{\operatorname{rk}}}
\nc{\Sets}{{{\cal{S}}ets}}
\nc{\Sew}{{{\cal{S}}ew}}
\nc{\sgn}{{\operatorname{sgn}}}
\nc{\Sh}{{{\cal S}h}}
\nc{\Sign}{{{\cal S}ign}}
\nc{\Spe}{{\mbox{\bf{Sp}}}}
\nc{\supr}{{\operatorname{sup}}}
\nc{\Supp}{{\operatorname{Supp}}}
\nc{\supp}{{\operatorname{supp}}}
\nc{\Teich}{{{\cal{T}}eich}}
\nc{\tFS}{{\widetilde{\cal{FS}}}}
\nc{\Tor}{{\operatorname{Tor}}}
\nc{\totimes}{{\tilde{\otimes}}}
\nc{\tr}{{\operatorname{tr}}}
\nc{\tRep}{{\widetilde{{\cal R}ep}}}
\nc{\tTeich}{{\widetilde{{\cal T}eich}}}
\nc{\Vc}{{\mbox{\bf{V}}_{\mbox{\bf{c}}}}}
\nc{\Vect}{{{\cal V}ect}}
\nc{\Ve}{{\operatorname{Vert}}}
\nc{\wt}{{\widetilde}}

\nc{\bo}{{\mbox{\bf{0}}}}
\nc{\One}{{\mbox{\bf{1}}}}
\nc{\one}{{\mbox{\bf{1}}}}

\nc{\BA}{{\Bbb A}}
\nc{\bA}{{\overline{A}}}
\nc{\ba}{{\mbox{\bf{a}}}}
\nc{\baB}{{\overline{B}}}
\nc{\baeta}{{\bar{\eta}}}
\nc{\baJ}{{\bar{J}}}
\nc{\BB}{{\Bbb B}}
\nc{\bB}{{\mbox{\bf{B}}}}
\nc{\bc}{{\mbox{\bf{c}}}}
\nc{\bC}{{\overline{C}}}
\nc{\BC}{{\Bbb{C}}}
\nc{\bCC}{{\overline{\cal{C}}}}
\nc{\bCM}{{\overline{\cal{M}}}}
\nc{\bD}{{\bar{D}}}
\nc{\BD}{{\overline{D}}}
\nc{\bd}{{\mbox{\bf{d}}}}
\nc{\BE}{{\overline{E}}}
\nc{\BF}{{\overline{F}}}
\nc{\bF}{{\mbox{\bf{F}}}}
\nc{\bg}{{\mbox{\bf{g}}}}
\nc{\bG}{{\mbox{\bf{G}}}}
\nc{\BG}{{\Bbb G}}
\nc{\bGamma}{{\overline{\Gamma}}}
\nc{\bbH}{{     {\mbox{\bf{H}}}_a       }}
\nc{\bH}{{\mbox{\bf{H}}}}
\nc{\bI}{{\mbox{\bf{I}}}}
\nc{\bL}{{\mbox{\bf{L}}}}
\nc{\BL}{{\Bbb{L}}}
\nc{\blambda}{{\bar{\lambda}}}
\nc{\bM}{{\mbox{\bf{V}}}}
\nc{\bmu}{{\vec{\mu}}}
\nc{\bN}{{\mbox{\bf{N}}}}
\nc{\BN}{{\Bbb{N}}}
\nc{\bnu}{{\mbox{\boldmath{${\nu}$}}}}
\nc{\bof}{{\mbox{\bf{f}}}}
\nc{\BP}{{\Bbb P}}
\nc{\bP}{{\mbox{\bf{P}}}}
\nc{\BPO}{{\overset{\circ}{\BP}}}
\nc{\BQ}{{\Bbb Q}}
\nc{\BR}{{\Bbb{R}}}
\nc{\bR}{{\mbox{\bf{R}}}}
\nc{\bp}{{\mbox{\bf{p}}}}
\nc{\barq}{{\mbox{\bf{q}}}}
\nc{\br}{{\mbox{\bf{r}}}}
\nc{\breta}{{\bar{\eta}}}
\nc{\bs}{{\mbox{\bf{s}}}}
\nc{\bS}{{\mbox{\bf{S}}}}
\nc{\bt}{{\mbox{\bf{t}}}}
\nc{\BT}{{\Bbb T}}
\nc{\bU}{{\mbox{\bf{U}}}}
\nc{\bV}{{\mbox{\bf{V}}}}
\nc{\bu}{{\mbox{\bf{u}}}}
\nc{\BUpsilon}{{\bar{\Upsilon}}}
\nc{\bw}{{\mbox{\bf{w}}}}
\nc{\bx}{{\mbox{\bf{x}}}}
\nc{\bX}{{\mbox{\bf{X}}}}
\nc{\BZ}{{\Bbb{Z}}}
\nc{\bz}{{\mbox{\bf{z}}}}
\nc{\bZ}{{\mbox{\bf{Z}}}}
\nc{\bzero}{\mbox{\boldmath{$0$}}}

\nc{\CA}{{\cal A}}
\nc{\CAD}{{\overset{\bullet}{\cal{A}}}}
\nc{\CAO}{{\overset{\circ}{\cal{A}}}}
\nc{\CB}{{\cal B}}
\nc{\CC}{{\cal C}}
\nc{\CalD}{{\cal D}}
\nc{\CE}{{\cal E}}
\nc{\CF}{{\cal F}}
\nc{\CG}{{\cal G}}
\nc{\CH}{{\cal H}}
\nc{\CI}{{\cal I}}
\nc{\CID}{{\overset{\bullet}{\cal{I}}}}
\nc{\CJ}{{\cal J}}
\nc{\CK}{{\cal K}}
\nc{\CL}{{\cal L}}
\nc{\CM}{{\cal M}}
\nc{\CN}{{\cal N}}
\nc{\CO}{{\cal O}}
\nc{\CP}{{\cal P}}
\nc{\CPO}{{\overset{\circ}{\cal{P}}}}
\nc{\CQ}{{\cal Q}}
\nc{\CR}{{\cal R}}
\nc{\CS}{{\cal S}}
\nc{\CT}{{\cal T}}
\nc{\CTD}{{\overset{\bullet}{\cal{T}}}}
\nc{\CTPO}{{\overset{\circ}{\cal{T}\cal{P}}}}
\nc{\CU}{{\cal{U}}}
\nc{\CV}{{\cal V}}
\nc{\CW}{{\cal W}}
\nc{\CX}{{\cal X}}
\nc{\CY}{{\cal Y}}
\nc{\CZ}{{\cal Z}}

\nc{\dCL}{{\overset{\bullet}{\cal{L}}}}
\nc{\dd}{{\operatorname{d}}}
\nc{\ddelta}{{\overset{\bullet}{\delta}}}
\nc{\dfu}{{\overset{\bullet}{\frak{u}}}}
\nc{\dlambda}{{\overset{\bullet}{\lambda}}}
\nc{\DO}{{\overset{\circ}{D}}}
\nc{\dpar}{{\partial}}
\nc{\dS}{{\overset{\bullet}{S}}}
\nc{\dT}{{\overset{\bullet}{T}}}

\nc{\fA}{{\frak{A}}}
\nc{\fb}{{\frak{b}}}
\nc{\fC}{{\frak{C}}}
\nc{\fD}{{\frak{D}}}
\nc{\fE}{{\frak{E}}}
\nc{\fF}{{\frak{F}}}
\nc{\fEF}{{\frak{EF}}}
\nc{\fFE}{{\frak{FE}}}
\nc{\ff}{{\frak{f}}}
\nc{\fg}{{\frak{g}}}
\nc{\fG}{{\frak{G}}}
\nc{\fH}{{\frak{H}}}
\nc{\fl}{{\frak{l}}}
\nc{\fK}{{\frak{K}}}
\nc{\fL}{{\frak{L}}}
\nc{\fM}{{\frak{M}}}
\nc{\fN}{{\frak{N}}}
\nc{\fn}{{\frak{n}}}
\nc{\fp}{{\frak{p}}}
\nc{\fu}{{\frak{u}}}
\nc{\fZ}{{\frak{Z}}}

\nc{\hCH}{{\hat{\cal{H}}}}
\nc{\hCI}{{\hat{\cal{I}}}}
\nc{\hfC}{{\hat{\frak{C}}}}
\nc{\hfg}{{\hat{\frak{g}}}}
\nc{\hL}{{\hat{L}}}
\nc{\HO}{{\overset{\circ}{H}}}
\nc{\hpsi}{{\hat{\psi}}}
\nc{\hx}{{\hat{x}}}

\nc{\jo}{{\overset{\circ}{j}}}

\nc{\phid}{{\overset{\bullet}{\phi}}}

\nc{\tA}{{\tilde{A}}}
\nc{\ta}{{\tilde{a}}}
\nc{\tB}{{\tilde{B}}}
\nc{\tb}{{\tilde{b}}}
\nc{\tBP}{{\tilde{\BP}}}
\nc{\tC}{{\tilde{C}}}
\nc{\tc}{{\tilde{c}}}
\nc{\tCA}{{\tilde{\cal{A}}}}
\nc{\tCC}{{\tilde{\cal{C}}}}
\nc{\tCH}{{\tilde{\cal{H}}}}
\nc{\tCI}{{\tilde{\cal{I}}}}
\nc{\tCO}{{\tilde{\cal{O}}}}
\nc{\tCP}{{\tilde{\cal{P}}}}
\nc{\tCT}{{\tilde{\cal{T}}}}
\nc{\tD}{{\tilde{D}}}
\nc{\tDelta}{{\tilde{\Delta}}}
\nc{\tE}{{\tilde E}}
\nc{\tF}{{\tilde F}}
\nc{\tfD}{{\tilde{\frak{D}}}}
\nc{\tfF}{{\tilde{\frak{F}}}}
\nc{\tff}{{\tilde{\frak{f}}}}
\nc{\tfu}{{\tilde{\frak{u}}}}
\nc{\tJ}{{\tilde{J}}}
\nc{\tj}{{\tilde{j}}}
\nc{\tK}{{\tilde K}}
\nc{\tL}{{\tilde{L}}}
\nc{\tM}{{\tilde{M}}}
\nc{\tP}{{\tilde{P}}}
\nc{\tPhi}{{\tilde{\Phi}}}
\nc{\tpi}{\tilde{\pi}}
\nc{\TPO}{{\overset{\circ}{T\BP}}}
\nc{\tR}{{\tilde{R}}}
\nc{\tS}{{\tilde S}}
\nc{\tT}{{\tilde{T}}}
\nc{\ttau}{{\tilde{\tau}}}
\nc{\ttheta}{{\tilde{\theta}}}
\nc{\tU}{{\tilde{U}}}
\nc{\tUpsilon}{{\tilde{\Upsilon}}}
\nc{\tW}{{\tilde W}}
\nc{\ty}{{\tilde y}}
\nc{\tY}{{\tilde Y}}
\nc{\txi}{{\tilde{\xi}}}

\nc{\UD}{{\overset{\bullet}{U}}}
\nc{\UO}{{\overset{\circ}{U}}}

\nc{\vA}{{\vec{A}}}
\nc{\valpha}{{\vec{\alpha}}}
\nc{\vbeta}{{\vec{\beta}}}
\nc{\vc}{{\vec{c}}}
\nc{\vD}{{\vec{D}}}
\nc{\vd}{{\vec{d}}}
\nc{\vgamma}{{\vec{\gamma}}}
\nc{\vK}{{\vec{K}}}
\nc{\vlambda}{{\vec{\lambda}}}
\nc{\vmu}{{\vec{\mu}}}
\nc{\vnu}{{\vec{\nu}}}
\nc{\vo}{{\vec{0}}}
\nc{\vu}{{\vec{u}}}
\nc{\vx}{{\vec{x}}}
\nc{\vy}{\vec{y}}
\nc{\vzero}{\vec{0}}

\nc{\XO}{{\overset{\circ}{X}}}

\nc{\ya}{{\operatorname{aj}}}

\nc{\nen}{\newenvironment}
\nc{\ol}{\overline}
\nc{\ul}{\underline}
\nc{\ra}{\rightarrow}
\nc{\lra}{\longrightarrow}
\nc{\Lra}{\Longrightarrow}
\nc{\lla}{\longleftarrow}
\nc{\Llra}{\Longleftrightarrow}
\nc{\hra}{\hookrightarrow}
\nc{\iso}{\overset{\sim}{\lra}}
\nc{\rlh}{\rightleftharpoons}

\nc{\IC}{{\cal{IC}}}
\nc{\PS}{{\cal{PS}}}
\nc{\oCQ}{{\overline{\cal Q}}}
\nc{\oCZ}{{\overline{\cal Z}}}
\nc{\dZ}{{\overset{\bullet}{\cal Z}}{}}
\nc{\oZ}{{\overset{\circ}{\cal Z}}{}}
\nc{\dP}{{\overset{\bullet}{\cal P}}{}}
\nc{\oP}{{\overset{\circ}{\cal P}}{}}
\nc{\Ue}{{U_\varepsilon}}
\nc{\Upe}{{\Upsilon_\varepsilon}}
\nc{\crho}{{\check{\rho}}}
\nc{\ctheta}{{\check{\theta}}}

\nc{\QL}[1]{\CQ^L_{#1}}
\nc{\QK}[1]{\CQ^K_{#1}}
\nc{\al}{\alpha}
\nc{\ga}{\gamma}
\nc{\ka}{\kappa}
\nc{\pr}{\bp\times \br}
\nc{\ME}{{\check \fE}}
\nc{\PP}{{\Bbb P}}
\nc{\CCC}{{\Bbb C}}
\nc{\NNN}{{\Bbb N}}
\nc{\ZZZ}{{\Bbb Z}}
\nc{\RT}{{\frak R\frak T}}
\nc{\Gr}{\operatorname{\text{Gr}}}
\nc{\wti}{\widetilde}
\nc{\wha}{\widehat}
\nc{\vphi}{\varphi}
\nc{\deff}{\operatorname{\text{def}}}
\nc{\hnu}{{\hat\nu}}
\nc{\tnu}{{\tilde\nu}}
\nc{\tka}{{\tilde\ka}}
\nc{\ti}{\tilde}
\nc{\fS}{{\frak S}}
\nc{\EGo}{{\overset{\circ}{\fE^\al_\Gamma}}}
\nc{\eac}{\fE^\al_{\underbrace{\theta_1,\dots,\theta_1}_{c_1},\dots,
\underbrace{\theta_\nu,\dots,\theta_\nu}_{c_\nu}}}
\nc{\lbr}{\{\!\{}
\nc{\rbr}{\}\!\}}

\nc{\Thm}[1]{Theorem~\ref{#1}}
\nc{\Prop}[1]{Proposition~\ref{#1}}
\nc{\Lem}[1]{Lemma~\ref{#1}}
\nc{\Cor}[1]{Corollary~\ref{#1}}
\nc{\Conj}[1]{Conjecture~\ref{#1}}
\nc{\Claim}[1]{Claim~\ref{#1}}
\nc{\Defn}[1]{Definition~\ref{#1}}
\nc{\Exa}[1]{Example~\ref{#1}}
\nc{\Rem}[1]{Remark~\ref{#1}}
\nc{\Note}[1]{Note~\ref{#1}}

\nen{thm}[1]{\label{#1}{\bf Theorem.\ } \em}{}
\nen{prop}[1]{\label{#1}{\bf Proposition.\ } \em}{}
\nen{lem}[1]{\label{#1}{\bf Lemma.\ } \em}{}
\nen{cor}[1]{\label{#1}{\bf Corollary.\ } \em}{}
\nen{conj}[1]{\label{#1}{\bf Conjecture.\ } \em}{}
\nen{claim}[1]{\label{#1}{\bf Claim.\ } \em}{}

\nen{defn}[1]{\label{#1}{\bf Definition.\ } }{}
\nen{exa}[1]{\label{#1}{\bf Example.\ } }{}

\nen{rem}[1]{\label{#1}{\em Remark.\ } }{}
\nen{note}[1]{\label{#1}{\em Note.\ } }{}
\nen{exer}[1]{\label{#1}{\em Exercise.\ } }{}

\begin{document}

\title[]{Global Intersection Cohomology of Quasimaps' Spaces}
\author{Michael Finkelberg}
\address{Independent Moscow University, 11 Bolshoj Vlasjevskij pereulok,
Moscow 121002 Russia}
\email{fnklberg@@main.mccme.rssi.ru}
\author{Alexander Kuznetsov}
\address{Independent Moscow University, 11 Bolshoj Vlasjevskij pereulok,
Moscow 121002 Russia}
\email{sasha@@ium.ips.ras.ru kuznetsov@@mpim-bonn.mpg.de}
\thanks{Both authors were partially supported by CRDF grant RM1-265.
M.F.~was partially supported by INTAS-94-4720}
\date{February 1997}
\maketitle

\section{Introduction}

\subsection{}
Let $C$ be a smooth projective curve of genus 0. Let $\CB$ be the variety
of complete flags in an $n$-dimensional vector space $V$.
Given an $(n-1)$-tuple $\alpha\in\BN[I]$
of positive integers one can consider the space $\CQ_\alpha$ of algebraic
maps of degree $\alpha$ from $C$ to $\CB$. This space is noncompact. Some
remarkable compactifications $\CQ^D_\alpha$ (Quasimaps),
$\CQ^L_\alpha$ (Quasiflags), $\CQ^K_\alpha$ (Stable Maps) of
$\CQ_\alpha$ were constructed by Drinfeld, Laumon and Kontsevich respectively.
In ~\cite{k} it was proved that the natural map $\pi:\ \CQ^L_\alpha\to
\CQ^D_\alpha$ is a small resolution of singularities. The aim of the present
note is to study the cohomology $H^\bullet(\CQ^L_\alpha,\BQ)$ of Laumon's
spaces or, equivalently, the Intersection Cohomology
$H^\bullet(\CQ^L_\alpha,IC)$ of Drinfeld's Quasimaps' spaces.

\subsection{}
It appears that $\CQ^L_\alpha$ admits a cell decomposition, whence its
cohomology has a pure Tate Hodge structure.
It was essentially computed by G.Laumon (see ~\cite{la}, Theorem 3.3.2).
For the reader's convenience we reproduce the computation in section 2.
We calculate the generating
function $P_G(t)$ (``Poincar\'e polynomial'') of the direct sum
$\oplus_{\alpha\in\BN[I]}H^\bullet(\CQ^D_\alpha,IC)$ as a formal cocharacter
of $G=SL_n$ with coefficients in the Laurent polynomials in $t$
(a formal variable of degree 2). It is given by the following formula:
$$P_G(t)=\frac{e^{2\rho}t^{-\frac{1}{2}\dim\CB}\sum_{w\in W}t^{\ell(w)}}
{\prod_{\theta\in R^+}(1-te^\theta)(1-t^{-1}e^\theta)}$$
where $W=S_n$ is the Weyl group of $G$ with its standard length function,
$R^+$ is the set of positive coroots of $G$, and
$2\rho$ stands for $\sum_{\theta\in R^+}\theta$.

\subsection{}
For any $\alpha,\gamma\in\BN[I]$ there is a closed subvariety of middle
dimension $\fE^\alpha_\gamma\subset\CQ^L_\alpha\times\CQ^L_{\alpha+\gamma}$.
It is formed by pairs of quasiflags such that the second one is a subflag
of the first one. The top-dimensional irreducible components of
$\fE^\alpha_\gamma$ are naturally numbered by the Kostant partitions
$\bc\in\fK(\gamma)$ of $\gamma$, independently of $\alpha$. For
$\bc\in\fK(\gamma)$ the corresponding irreducible component $\fE^\alpha_\bc$,
viewed as a correspondence between $\CQ^L_\alpha$ and $\CQ^L_{\alpha+\gamma}$,
defines two operators:
$$e_\bc:\ H^\bullet(\CQ^L_\alpha)\rightleftharpoons
H^\bullet(\CQ^L_{\alpha+\gamma})\ :f_\bc$$
adjoint to each other with respect to Poincar\'e duality.

\subsubsection{}
Let $\fn$ denote the Lie subalgebra of upper-triangular matrices in the Lie
algebra ${\frak{sl}}_n$. Let $U(\fn)$ denote the Kostant integral form
(with divided powers) of the universal enveloping algebra of $\fn$. It turns
out that the linear span of operators $e_\bc$ is closed under composition;
the algebra they form is naturally isomorphic to $U(\fn)$, and the isomorphism
takes $e_\bc$ to the corresponding element of Poincar\'e-Birkhoff-Witt-Kostant
basis of $U(\fn)$.

\subsubsection{}
Moreover, all the operators $e_\bc,f_\bc$ together generate an action of
the universal enveloping algebra of ${\frak{sl}}_n$ on
$\oplus_{\alpha\in\BN[I]}H^\bullet(\CQ^L_\alpha,\BQ)$. The character of
this module is given by
$\dfrac{|W|e^{2\rho}}{\prod_{\theta\in R^+}(1-e^\theta)^2}$.

\subsection{}
We conjecture that the ${\frak{sl}}_n$-module
$\oplus_{\alpha\in\BN[I]}H^\bullet(\CQ^L_\alpha,\BQ)$ is isomorphic to
$H^\nu_\fn(\CN,\CO)$. Here $\CN$ stands for the nilpotent cone of
${\frak{sl}}_n$, and $H^\nu_\fn(\CN,\CO)$ is the cohomology of the structure
sheaf with supports in $\fn$ of degree $\nu=\dim\fn=\frac{n(n-1)}{2}$.
To verify this conjecture it would be enough to check that
$\oplus_{\alpha\in\BN[I]}H^\bullet(\CQ^L_\alpha,\BQ)$ is free as a
$U(\fn)$-module (see section 6).

\subsection{}
The above conjecture is motivated by B.Feigin's conjecture about the
semiinfinite cohomology $H_\fu^{\frac{\infty}{2}+\bullet}$
of small quantum group $\fu$ of type $A_{n-1}$ (see ~\cite{ar} and section 6).
Let us add a few more words about motivation.

We believe that the Drinfeld's spaces $\CQ^D_\alpha$ are the basic building
blocks of the would-be {\em Semiinfinite Flag Space} (cf. ~\cite{fm}).
On the other hand, it was conjectured by G.Lusztig and B.Feigin that an
appropriate category of perverse sheaves on Semiinfinite Flags is equivalent
to a regular block of the category $\CC$ of graded $\fu$-modules. In this
equivalence, the algebraic counterpart of the
global Intersection Cohomology $H^\bullet(\CQ^D_\alpha,IC)$ is exactly
$_{(\alpha+2\rho)}H_\fu^{\frac{\infty}{2}+\bullet}$ ---
the (co)weight $(\alpha+2\rho)$
space of the ${\frak{sl}}_n$-module of semiinfinite
cohomology $H_\fu^{\frac{\infty}{2}+\bullet}$.

In fact, another geometric realization of the category $\CC$ was constructed
in ~\cite{fs}. One of the main theorems of {\em loc. cit.}
canonically identifies $_{(\alpha+2\rho)}H^{\frac{\infty}{2}+\bullet}_\fu$
with the Intersection Cohomology of a certain one-dimensional local system
on a certain configuration space of $C$.

Combining all the established equalities of characters (see section 6)
we see that the above Intersection Cohomology has the
same graded dimension as $H^\bullet(\CQ^D_\alpha,IC)$.

We believe that it would be extremely interesting and important to find
a direct explanation of this coincidence. In fact, this (conjectural)
coincidence was the main impetus for the present work. The desired
explanation might be not that easy since
$H^\bullet(\CQ^D_\alpha,IC)$ has a Tate Hodge structure while the
Intersection Cohomology of the above local system has quite a nontrivial
Hodge structure (e.g. of elliptic curves or K3-surfaces) already in the
simplest examples.

\subsection{}
The idea to realize the algebra $U(\fn)$ in correspondences (or in the
$K$-groups of constructible sheaves on certain spaces) is not new:
see e.g. the remarkable works ~\cite{lu}, ~\cite{blm}, ~\cite{gi}, ~\cite{naa}.
What seems to be new compared to {\em loc. cit.} is the reason behind the
relations in $U(\fn)$ (or $U({\frak{sl}}_n)$). Say, the divided powers of
simple generators appear in {\em loc. cit.}, roughly, due to the fact that
the flag manifold of $SL_d$ has Euler characteristic $d!$; while in the
present work the divided powers appear, roughly, due to the fact that
the Cartesian power $C^d$ is a $d!$-fold cover of the symmetric power
$C^{(d)}$.

Thus, the present construction may be viewed as a sort of globalization of
{\em loc. cit.} in the particular case of Dynkin diagram of type $A_{n-1}$
(one might say that ~\cite{blm} and ~\cite{gi} lived in the formal
neighbourhood of a point $0\in C$, while we work over the whole $C$).
Note that the constructions of ~\cite{lu} and ~\cite{naa}
can be (and are) generalized to arbitrary Dynkin graphs and quantized. It would
be extremely interesting to generalize the results of the present note to an
arbitrary Dynking graph (or even quantize them).

On the other hand, the idea to realize Lie algebras' representations via
``global'' correspondences is not new either: see e.g. the remarkable
works ~\cite{gr}, ~\cite{na}. These works
realize some irreducible representations of infinite
dimensional Lie algebras (Heisenberg and Clifford) in the cohomology of
Hilbert schemes of {\em surfaces}. Thus the present work may be viewed
as a baby version of {\em loc. cit.} Note though that in all the previous cases
the representations of Lie algebras realized geometrically turned out to be
irreducible, while we expect our modules to be {\em nonsemisimple}. In fact,
we conjecture that they are {\em tilting} (see ~\cite{ap} and section 6).
Also, in the global context, the appearence of Serre relations seems to be new.

\subsection{}
It is clear from the above explanations how much we were influenced by all
the above cited works. It is a pleasure to thank B.Feigin, S.Arkhipov,
and V.Ostrik for the numerous illuminating discussions and suggestions.
Above all we are obliged to I.Mirkovic who spent half a year teaching one
of us (M.F.) the beautiful geometry of affine and semiinfinite flag spaces,
of which the present results are but a superficial manifestation.

\subsection{}
The present note is a sequel to ~\cite{k}. We will freely refer the reader
to {\em loc. cit.}

\section{Cohomology of $\CQ^L_\alpha$}

\subsection{Notations}

\subsubsection{}
\label{not}
We choose a basis $\{v_1,\ldots,v_n\}$ in $V$. This choice
defines a Cartan subgroup $H\subset G$ of matrices diagonal with respect to
this basis, and a Borel subgroup $B\subset G$ of matrices upper triangular
with respect to this basis. We have $\CB=G/B$.

Let $I=\{1,\ldots,n-1\}$ be the set of simple coroots of $G=SL_n$.
Let $R^+$ denote the set of positive coroots,
and let $2\rho=\sum_{\theta\in R^+}\theta$.
For $\alpha=\sum a_ii\in\BN[I]$ we set $|\alpha|:=\sum a_i$.

Recall the notations of ~\cite{k} concerning Kostant's partition function.
For $\gamma\in\BN[I]$ a {\em Kostant partition} of $\gamma$ is a decomposition
of $\gamma$ into a sum of positive coroots with multiplicities.
The set of Kostant partitions of $\gamma$ is denoted by
$\fK(\gamma)$. For $\kappa\in\fK(\gamma)$ let $|\kappa|=\gamma$,
$||\kappa||=|\gamma|$ and let $K(\kappa)$ be the number of summands in $\kappa$.

There is a natural bijection between the set of pairs $1\leq q\leq p\leq n-1$
and $R^+$, namely, $(p,q)$ corresponds to $i_q+i_{q+1}+\ldots+i_p$. Thus a
Kostant partition $\kappa$ is given by a collection of nonnegative integers
$(\kappa_{p,q}), 1\leq q\leq p\leq n-1$.
Following {\em loc. cit.} (9) we define a collection $\mu(\kappa)$ as follows:
$\mu_{p,q}=\sum_{r\leq q\leq p\leq s}\kappa_{s,r}$.

\subsubsection{}
For the definition of Laumon's Quasiflags' space $\CQ^L_\alpha$ the reader
may consult ~\cite{la} 4.2, or ~\cite{k} 1.4. It is the space of complete
flags of locally free subsheaves
$$0\subset E_1\subset\dots\subset E_{n-1}\subset V\otimes\CO_C$$
such that rank$(E_k)=k$, and $\deg(E_k)=-a_k$.

It is known to be a smooth
projective variety of dimension $2|\alpha|+\dim\CB$.

\subsubsection{}
For the definition of Drinfeld's Quasimaps' space $\CQ^D_\alpha$ the
reader may consult ~\cite{k} 1.2. It is the space of collections of
invertible subsheaves $\CL_\lambda\subset V_\lambda\otimes\CO_C$ for
each dominant weight $\lambda\in X^+$ satisfying Pl\"ucker relations,
and such that $\deg\CL_\lambda=-\langle\lambda,\alpha\rangle$.

It is known to be a (singular, in general) projective variety of
dimension $2|\alpha|+\dim\CB$.

\subsection{}
Given a quasiflag $E_\bullet\in\CQ^L_\alpha$ and a point $x\in C$ the
{\em type $\kappa(E_\bullet),\mu(E_\bullet)$ (of defect) of $E$ at $x$}
was defined in {\em loc. cit.} (6)--(11).

{\bf Definition.} For $\gamma\leq\alpha,\ \kappa\in\fK(\gamma)$ we denote
by $\CZ^\kappa_\alpha\subset\CQ^L_\alpha$ the locally closed subspace formed
by the quasiflags with defect of type $\kappa$ at $\infty\in C$. In particular,
$\CZ^0_\alpha$ is an open subset of $\CQ^L_\alpha$.

{\em Normalization} at $\infty\in C$ (see {\em loc. cit.} 1.5.1) defines
a map $$\varpi^\kappa_\alpha:\ \CZ^\kappa_\alpha\lra\CZ^0_{\alpha-|\ka|}$$

Evaluation at $\infty\in C$ defines a map
$$\Upsilon_\alpha:\ \CZ^0_\alpha\lra\CB$$
Evidently, $\Upsilon_\alpha$ is a locally trivial fibration. We will denote
the fiber of $\Upsilon_\alpha$ over the point $B\in\CB=G/B$ by $\CY_\alpha$.

\subsection{}
\label{retract}
The point $B\in\CB$ is represented by a flag $0\subset V_1\subset\ldots
\subset V_{n-1}\subset V$. Let $E^0_\bullet$ denote the corresponding
trivial flag of subbundles: $E^0_i=V_i\otimes\CO_C$.
For the point $0\in C$ the {\em simple fiber} $F(E^0_\bullet,\alpha0)\subset
\CQ^L_\alpha$ was introduced in {\em loc. cit.} 2.1.2. Its cohomology was
computed in {\em loc. cit.} 2.4.4: its Poincar\'e polynomial equals
$\CK_\alpha(t)=t^{|\alpha|}\sum_{\kappa\in\fK(\alpha)}t^{-K(\kappa)}$ ---
the Lusztig-Kostant polynomial ($t$ has degree 2).

{\bf Lemma.} The closed embedding $F(E^0_\bullet,\alpha0)\hookrightarrow
\CY_\alpha$ induces an isomorphism of cohomology:
$$H^\bullet(\CY_\alpha,\BQ)\iso H^\bullet(F(E^0_\bullet,\alpha0),\BQ)$$

{\em Proof.} We restrict the natural map $\pi:\ \CQ^L_\alpha\lra\CQ^D_\alpha$
to the locally closed subvariety $\CY_\alpha\subset\CQ^L_\alpha$.
The image $\pi(\CY_\alpha)\subset\CQ^D_\alpha$ is denoted by $\CZ_\alpha$.
It consists of quasimaps regular at $\infty\in C$ and taking there the
value $B\in\CB$. We will preserve the same name for the restriction of
$\pi$ to $\CY_\alpha\lra\CZ_\alpha$. It follows from the main Theorem of
{\em loc. cit.} that $\pi$ is a small resolution of singularities.
Hence $H^\bullet(\CY_\alpha,\BQ)=H^\bullet(\CZ_\alpha,IC)$.

Now $\BC^*$ acts on $C$ by dilations preserving $0,\infty$, and thus it
acts on both $\CY_\alpha$ and $\CZ_\alpha$, and the map $\pi$ is equivariant
with respect to this action. The space $\CZ_\alpha$ has the only point $Z$
fixed by $\BC^*$: it is the point where all the defect is concentrated at
$0\in C$. The {\em simple fiber} $F(E^0_\bullet,\alpha0)$ is nothing
else than the fiber $\pi^{-1}(Z)$. So the stalk $IC_{(Z)}$
 of $IC$-sheaf at the point
$Z$ equals $H^\bullet(F(E^0_\bullet,\alpha0),\BQ)$. On the other hand,
since the $\BC^*$-action contracts $\CZ_\alpha$ to $Z$, we have
$H^\bullet(\CZ_\alpha,IC)=IC_{(Z)}$. The proposition is proved.

Alternatively, instead of using Intersection Cohomology, we could argue
that according to ~\cite{sl}, 4.3, $F(E^0_\bullet,\alpha0)=\pi^{-1}(Z)$
is a deformation retract of $\CY_\alpha$. $\Box$

\subsubsection{Remark} The space $\CZ_\alpha$ plays a central role and is
extensively studied in ~\cite{fm}.

\subsubsection{Corollary}
\label{van}
The odd-dimensional cohomology of $\CY_\alpha$
vanishes.

{\em Proof.} Follows immediately from ~\cite{k}, 2.4.4. $\Box$

\subsection{}
\label{fib}
We consider the locally trivial fibration $\Upsilon_\alpha:\ \CZ^0_\alpha
\lra\CB$ with the fiber $\CY_\alpha$. Since the odd-dimensional cohomology
of both the fiber and the base vanishes, the Leray spectral sequence of
this fibration degenerates, and we arrive at the following Lemma:

{\bf Lemma.} The odd-dimensional cohomology of $\CZ^0_\alpha$ vanishes.
The Poincar\'e polynomial $P(H^\bullet(\CZ^0_\alpha),t)$ equals
$\CK_\alpha(t)\sum_{w\in W}t^{\ell(w)}$. $\Box$

\subsection{Lemma}
\label{fig}
The Poincar\'e polynomial of the cohomology with compact support
$P(H_c^\bullet(\CZ^0_\alpha),t)$ equals $t^{\dim\CB+2|\alpha|}
\CK_\alpha(t^{-1})\sum_{w\in W}t^{-\ell(w)}$.

{\em Proof.} The space $\CZ^0_\alpha$ is smooth of dimension
$\dim\CB+2|\alpha|$. Now apply the Poincar\'e duality and the Lemma ~\ref{fib}.
$\Box$

\subsection{Lemma}
\label{mig}
The odd-dimensional cohomology with compact support of $\CZ^0_\alpha$ vanishes.
The Poincar\'e polynomial of the cohomology with compact support
$P(H_c^\bullet(\CZ^\kappa_\alpha),t)$ equals
$$
t^{\dim\CB+2|\alpha|-||\kappa||-K(\kappa)}
\CK_{\alpha-|\kappa|}(t^{-1})\sum_{w\in W}t^{-\ell(w)}.
$$

{\em Proof.} The normalization map $\varpi^\kappa_\alpha:\
\CZ^\kappa_\alpha\lra\CZ^0_{\alpha-|\kappa|}$ is a locally trivial fibration
with a fiber isomorphic to a pseudoaffine space $\fS_{\mu(\kappa)}$
(see ~\cite{k}~(16)) of dimension $||\kappa||-K(\kappa)$
(see {\em loc. cit.}~(21)). Now apply the Lemma ~\ref{fig}. $\Box$

\subsection{}
\label{character}
We consider the stratification
$$\CQ^L_\alpha=\bigsqcup_{|\kappa|\leq\alpha}\CZ^\kappa_\alpha$$
and the corresponding Cousin spectral sequence converging to the compactly
supported cohomology of $\CQ^L_\alpha$ (equal to $H^\bullet(\CQ^L_\alpha,\BQ)$
by Poincar\'e duality). Since the odd-dimensional compactly supported cohomology
of every stratum vanishes, the Cousin spectral sequence degenerates, and we
are able to compute the Poincar\'e polynomial of the space
$\CQ^L_\alpha$. To write it down in a neat form we will need some minor
preparations.

First of all, we shift the cohomological degree so that the cohomology becomes
symmetric around zero degree: we consider
$H^\bullet(\CQ^L_\alpha,\BQ[\dim\CQ^L_\alpha])$. Recall that $\dim\CQ^L_\alpha=
2|\alpha|+\dim\CB=2|\alpha|+\frac{n(n-1)}{2}$.

Second, we will consider the generating function for
$\oplus_{\alpha\in\BN[I]}H^\bullet(\CQ^L_\alpha,\BQ[\dim\CQ^L_\alpha])$.
To record the information on $\alpha$ we will consider this generating
function as a formal cocharacter of $H$ with coefficients in the Laurent
polynomials in $t$. Formal cocharacters will be written multiplicatively,
so that the cocharacter corresponding to $\alpha$ will be denoted by $e^\alpha$.
Finally, for the reasons which will become clear later (see Proposition
~\ref{h}), we will make the following rescaling. We will attach to
$H^\bullet(\CQ^L_\alpha,\BQ[\dim\CQ^L_\alpha])$ the cocharacter
$e^{\alpha+2\rho}$.

With all this in mind, the Poincar\'e polynomial $P_G(t)$ of
$\oplus_{\alpha\in\BN[I]}H^\bullet(\CQ^L_\alpha,\BQ[\dim\CQ^L_\alpha])$
is calculated as follows:

{\bf Theorem.} $$P_G(t)=\frac{e^{2\rho}t^{-\frac{1}{2}\dim\CB}
\sum_{w\in W}t^{\ell(w)}}
{\prod_{\theta\in R^+}(1-te^\theta)(1-t^{-1}e^\theta)}$$
$\Box$

\subsection{}
\label{IC}
The main Theorem of ~\cite{k} asserts that the natural map
$\pi:\ \CQ^L_\alpha\lra\CQ^D_\alpha$ is a small resolution
of singularities. Hence the Intersection Cohomology complex on
$\CQ^D_\alpha$ is the direct image of the constant sheaf on $\CQ^L_\alpha:\
IC=\pi_*\underline{\BQ}[\dim\CQ^L_\alpha]$. This implies that the global
Intersection Cohomology $H^\bullet(\CQ^D_\alpha,IC)$ coincides with
the cohomology $H^\bullet(\CQ^L_\alpha,\BQ[\dim\CQ^L_\alpha])$.
Thus we obtain the following theorem.

{\bf Theorem.}
The generating
function for the global Intersection Cohomology of Drinfeld's Quasimaps'
spaces $\oplus_{\alpha\in\BN[I]}H^\bullet(\CQ^D_\alpha,IC)$ is given by
$$P_G(t)=\frac{e^{2\rho}t^{-\frac{1}{2}\dim\CB}\sum_{w\in W}t^{\ell(w)}}
{\prod_{\theta\in R^+}(1-te^\theta)(1-t^{-1}e^\theta)}$$
$\Box$

\subsection{}
\label{main}
{\bf Theorem.} There is a cell decomposition
$$\CQ^L_\alpha=\bigsqcup\Delta(w,\kappa^0,\kappa^\infty)$$
into cells numbered by the following data:
$w\in W$; partition $\kappa^0$ (resp.
$\kappa^\infty$) of $\gamma^0\in\BN[I]$ (resp. $\gamma^\infty\in\BN[I]$)
such that $\gamma^0+\gamma^\infty=\alpha$.

\subsubsection{} The proof of the Theorem will occupy the rest of the section.

\subsection{} We will consider a torus action on $\CQ^L_\alpha$ with finitely
many fixed points, and the Bialynicki-Birula decomposition defined by this
action will give the desired cell decomposition.

\subsubsection{}
\label{zvezdochka}
The Cartan group $H$ acts on $V$ and hence on $\CQ^L_\alpha$.
The group $\BC^*$ of dilations of $C={\Bbb P}^1$ preserving $0$ and $\infty$
also acts on $\CQ^L_\alpha$ commuting with the action of $H$.
Hence we obtain the action of a torus $\BT:=H\times\BC^*$ on $\CQ^L_\alpha$.

\subsubsection{} We fix a coordinate $z$ on $C={\Bbb P}^1$ such that
$z(0)=0,\ z(\infty)=\infty$.

\subsection{}
\label{fixpoint}
Let us describe the fixed point set $(\CQ^L_\alpha)^\BT$.

Given a triple $(w,\kappa^0,\kappa^\infty)$ as in the Theorem ~\ref{main},
we define the point $\delta(w,\kappa^0,\kappa^\infty)\in(\CQ^L_\alpha)^\BT$
as follows.

It is a quasiflag $(E_1,\ldots,E_{n-1})$ such that its
normalization $(\tilde{E}_1,\ldots,\tilde{E}_{n-1})$ (see ~\cite{k},
Definition ~1.5.1) is a constant flag with $\tilde{E}_1$ spanned by $v_{w(1)}$;
$\tilde{E}_2$ spanned by $v_{w(1)}$ and $v_{w(2)};\ \ldots;\ \tilde{E}_{n-1}$
spanned by $v_{w(1)},v_{w(2)},\ldots,v_{w(n-1)}$.

Its defect is a collection of torsion sheaves (see {\em loc. cit.}) on $C$
supported at $0$ and $\infty$.

In a neighbourhood of $0\in C$ the quasiflag $(E_1,\ldots,E_{n-1})$ is
defined as follows:

$$
{\arraycolsep=1pt
\begin{array}{llrlcrlcccrlc}
E_1 & =\langle & z^{d^0_{1,1}} & v_{w(1)}&\rangle \\
E_2 & =\langle & z^{d^0_{2,1}} & v_{w(1)}&,& z^{d^0_{2,2}} & v_{w(2)}
&\rangle\\
\ \vdots && \vdots &&& \vdots \\
E_{n-1} & =\langle & z^{d^0_{n-1,1}} & v_{w(1)} &,& z^{d^0_{n-1,2}} & v_{w(2)} 
& , & \dots & , & z^{d^0_{n-1,n-1}} & v_{w(n-1)}&\rangle \\
\end{array}
}
$$
where the collection
$(d^0_{p,q})_{1\leq q\leq p\leq n-1}$ (resp.
$(d^\infty_{p,q})_{1\leq q\leq p\leq n-1}$) is defined via $\kappa^0$
(resp. $\kappa^\infty$) as follows:
$$
d^\bullet_{p,q}=\sum_{r=p}^{n-1}\kappa^\bullet_{r,q}.
$$

Finally, in a neighbourhood of $\infty\in C$ the quasiflag
$(E_1,\ldots,E_{n-1})$ is defined exactly as around $0$, with the replacement
of $d^0$ by $d^\infty$ and $z$ by $z^{-1}$.

\subsection{Proposition} The fixed point set $(\CQ^L_\alpha)^\BT$
coincides with
the collection of points $\delta(w,\kappa^0,\kappa^\infty)$ numbered by the
triples as in the Theorem ~\ref{main}.

{\em Proof.} Easy. $\Box$

\subsection{}
\label{cells}
It is well known that for a generic choice of one-parametric subgroup
$\fG\subset \BT$ we have $(\CQ^L_\alpha)^\fG=(\CQ^L_\alpha)^\BT$. Hence we can
apply the main Theorem 4.4 of ~\cite{bb} to obtain the desired decomposition
$$\CQ^L_\alpha=\bigsqcup\Delta(w,\kappa^0,\kappa^\infty)$$ into locally
closed subchemes. Each one of them is an affine space containing exactly
one fixed point: namely, $\Delta(w,\kappa^0,\kappa^\infty)\ni
\delta(w,\kappa^0,\kappa^\infty)$. This completes the proof of the Theorem
~\ref{main}. $\Box$

\subsubsection{Remark}
It would be desirable to make a wise canonical choice of $\fG$ producing
some canonical cell decomposition (or, moreover, a stratification) of
$\CQ^L_\alpha$. For instance, we expect that for such a wise choice
the dimension of $\Delta(w,\kappa^0,\kappa^\infty)$ would be given by
$d(w,\kappa^0,\kappa^\infty)=\ell(w)+||\kappa^0||+||\kappa^\infty||+
K(\kappa^0)-K(\kappa^\infty)$. This would give a more natural proof of
the Theorem ~\ref{character}.

Furthermore, such decomposition would produce a canonical basis in
$H^\bullet(\CQ^L_\alpha,\BQ)$
(Poincar\'e duals of fundamental classes of cells)
which in turn might prove useful in checking the freeness of
$U(\fn)$-action on $\oplus_{\alpha\in\BN[I]}H^\bullet(\CQ^L_\alpha,\BQ)$
defined in the next section (see the Conjecture ~\ref{conjecture} and
the Remark ~\ref{freeness}).

Unfortunately, we were not able to make such a wise choice of $\fG$.

\section{Simple correspondences}

\subsection{}\label{def}
For any $i\in I$ and $\alpha\in\BN[I]$ we introduce the following closed
subvariety $\fE_i^\alpha\subset\CQ^L_\alpha\times\CQ^L_{\alpha+i}$.

{\bf Definition.} $\fE_i^\alpha:=
\{((E_1,\ldots,E_{n-1}),(E'_1,\ldots,E'_{n-1}))$
such that for $j\not=i$ we have $E_j=E'_j$, while $E_i\supset E'_i$, and
$E_i/E'_i$ is a torsion sheaf of length one$\}$.

There are natural maps
$$
\bp:\fE_i^\al\to\QL\al,\qquad
\barq:\fE_i^\al\to\QL{\al+i},\quad\text{and}\quad
\br:\fE_i^\al\to C.
$$
The first and second maps are induced by the projections of
$\QL\al\times\QL{\al+i}$ onto the first and second factor and the third
is defined as
$$
\br((E_\bullet,E'_\bullet))=\supp{E_i/E'_i}.
$$

The following Lemma describes the fibers of the map
$$
\bp\times \br:\fE_i^\al\to\QL\al.
$$

\subsection{Lemma}
\label{fibersofpr}
Let $E_\bullet\in\QL\al$, $x\in C$.
The fiber $(\bp\times \br)^{-1}(E_\bullet)$ is
naturally isomorphic to the projective space $\PP(\Hom(E_i/E_{i-1},\CO_x))$.
The map $(\bp\times \br)$ is an isomorphism over the space of
pairs $(E_\bullet,x)$ such that $E_\bullet\in\CQ_\al$ is
a flag of subbundles.

{\em Proof.}
The fiber $(\bp\times \br)^{-1}(E_\bullet,x)$
is evidently isomorphic to the space of all subsheaves
$E'_i\subset E_i$ such that $E_i/E'_i=\CO_x$ and embedding
$E_{i-1}\hookrightarrow E_i$ factors through $E'_i$.
In other words it is the space of all diagrams
$$
\begin{CD}
E_{i-1} @>>>    E'_i    @>>>    E'_i/E_{i-1}    \\
@|              @VVV            @VVV            \\
E_{i-1} @>>>    E_i     @>>>    E_i/E_{i-1}     \\
@.              @VVV            @VVV            \\
        @.      \CO_x   @=      \CO_x
\end{CD}
$$
with exact rows and columns. But such a diagram can be uniquely
reconstructed from the map $E_i/E_{i-1}\to\CO_x$.
The first part of the Lemma follows.

If $E_\bullet\in\CQ_\al$ then the quotient $E_i/E_{i-1}$
is locally free, hence $\Hom(E_i/E_{i-1},\CO_x)=\CCC$ for all
$x$. This means that $\bp\times \br$ is an isomorphism over
the space $\CQ_\al\times C\subset\QL\al\times C$. $\Box$

\subsubsection{} Let $\fE_{\lbr i\}\}}^\al$ be the closure of the
space $(\bp\times \br)^{-1}(\CQ_\al\times C)$.
Recall (see ~\cite{k}, 1.4.1) that $\dim\CQ^L_\alpha=\dim\CB+2|\alpha|$.
The map $(\bp\times \br)_{|\fE_{\lbr i\}\}}^\al}:
\fE_{\lbr i\}\}}^\al\to\QL\al\times C$
is birational, hence $\fE_{\lbr i\}\}}^\al$ is a
$(\dim\CB+2|\al|+1)$-dimensional irreducible variety.

\begin{lem}{irred}
The space $\fE_{\lbr i\rbr }^\al$ is a unique $(\dim\CB+2|\al|+1)$-dimensional
irreducible component of $\fE_i^\al$.
\end{lem}

{\em Proof.}
It is a particular case of the Proposition ~\ref{dims}.

Alternatively, a direct proof goes as follows.

Consider the following stratification of $\QL\al\times C$
$$
\QL\al\times C=\bigsqcup\begin{Sb}\gamma\le\alpha\\\ka\in\fK(\gamma)\end{Sb}
Z^\al_\ka,
$$
where $Z^\al_\ka\subset\QL\al\times C$ is the subspace of pairs
$(E_\bullet,x)$ such that $\ka\in\fK(\ga)$ is the type of the defect of
$E_\bullet\in\QL\al$ at the point $x\in C$ (see the section 2 of ~\cite{k}).
Considering the map
$$
\Pi:Z^\al_\ka\to\QL{\al-|\ka|}\times C,\qquad (E_\bullet,x)\mapsto
(\ti E_\bullet,x),
$$
where $\ti E_\bullet$ is the normalization at $x$ of $E_\bullet$
and applying the Lemma~2.4.3
of {\em loc.\ cit.} we see that
$$
\dim Z^\al_\ka=\dim\CB+2|\al|-||\ka||-K(\ka)+1.
$$
On the other hand, it is easy to see that over the stratum $Z^\al_\ka$
we have
$$
\length(T_x)=\sum_{p=1}^{i-1}\ka_{i-1,p},
$$
where $T_x$ is the part of the torsion $T$ in the quotient sheaf
$E_i/E_{i-1}$ with support at $x$. Using obvious inequality
$$
\dim\Hom(T,\CO_x)\le\length(T_x)
$$
we see that the dimension of the fiber of $\fE_i^\al$ over the
point $(E_\bullet,x)\in Z^\al_\ka$ is not greater than
$\sum_{p=1}^{i-1}\ka_{i-1,p}$. But
$$
||\ka||+K(\ka)\ge\sum_{p=1}^{i-1}\ka_{i-1,p}
$$
and equality is possible only for $\gamma=0$ (it follows easily from
{\em loc.\ cit.}~(9)). This means that for any $0<\gamma\le\alpha$
and $\ka\in\fK(\gamma)$
$$
\dim \bp^{-1}(Z^\al_\ka)<\dim\CB+2|\al|+1
$$
and the Lemma follows.  $\Box$

\subsection{}
\label{ef}
According to Lemma ~\ref{irred}, we may consider the Poincar\'e dual
of the fundamental class
$[\fE_{\lbr i\rbr }^\alpha]\in H^\bullet(\CQ^L_\alpha\times
\CQ^L_{\alpha+i},\BQ)$. Viewed as a correspondence, it defines two operators:
$$e_i:\ H^\bullet(\CQ^L_\alpha,\BQ)\rightleftharpoons
H^\bullet(\CQ^L_{\alpha+i},\BQ)\ :f_i$$
adjoint to each other with respect to Poincar\'e duality.
The operator $e_i$ increases the cohomological degree by 2,
and the operator $f_i$ decreases it by 2.

\subsubsection{Remark}\label{rem1}
We may also consider the operators 
$$
\hat e_i:\ H^\bullet(\CQ^L_\alpha,\BQ)\rightleftharpoons
H^\bullet(\CQ^L_{\alpha+i},\BQ)\ :\hat f_i
$$
defined by the fundamental class 
$[\fE_i^\al]\in H^\bullet(\CQ^L_\alpha\times\CQ^L_{\alpha+i},\BQ)$.
Consider the decompositions of the operators $\hat e_i$ and $\hat f_i$ into
the sum of operators shifting the cohomological degree by $k$.
$$
\hat e_i=\sum\hat e_i^k,\qquad
\hat f_i=\sum\hat f_i^k
$$
The Lemma~\ref{irred} implies that
$$
e_i=\hat e_i^2,\qquad
f_i=\hat f_i^{-2},
$$
and
$$
\hat e_i^k=0\quad\text{for }k>2\qquad\text{and}\qquad
\hat f_i^k=0\quad\text{for }k<-2.
$$

\subsubsection{}
We fix an orientation $\Omega=(1\lra2\lra\ldots\lra n-1)$
of the Dynkin graph with the set of vertices $I$.
Note that any flag of subsheaves (subbundles) in the trivial bundle 
$V\otimes\CO_C$ can be considered as a representation of the quiver
$\Omega$ in the category of subsheaves (subbundles) of $V\otimes\CO_C$.
Therefore, given a pair of flags $E'_\bullet\subset E_\bullet$
we have the quotient representation $T_\bullet=E_\bullet/E'_\bullet$.
This is a representation of the quiver $\Omega$ in the category
of torsion sheaves on $C$. Let us denote the category of
such representations by $\RT$. Define the {\em dimension} and {\em local
dimension at $x\in C$} of
$T=(T_1,\dots,T_{n-1})\in\Ob(\RT)$ as the coroots
$$
\dim T=\sum_{i\in I}\length(T_i)i\in\NNN[I],\qquad
\dim_x T=\sum_{i\in I}\length_x(T_i)i\in\NNN[I].
$$
Let $T\in\Ob(\RT)$ and $\dim T=\gamma$.
Given a filtration $0=F_0\subset F_1\dots\subset F_m=T$ of $T$ by
subrepresentations we say that it is a filtration of the type 
$(\gamma_1,\dots,\gamma_m)$ if $\dim F_k/F_{k-1}=\gamma_k$.

Let us denote by $\CO_x[i]$ a simple $i$-dimensional object
in the category $\RT$, consisting of the sheaf $\CO_x$
which lives over the $i$-th vertex of $\Omega$.

\subsection{Proposition}
\label{triv}
Given $i,j\in I$ such that $|i-j|>1$ we have

a) $e_ie_j=e_je_i$;

b) $f_if_j=f_jf_i$.

{\em Proof.} 
Instead of $e_ie_j$ and $e_je_i$ we will consider the
components of $\hat e_i\hat e_j$ and $\hat e_j\hat e_i$
shifting the cohomological degree by 4 (it suffices by the
Remark~\ref{rem1}). To this end, consider
the spaces
$$
\fE^\al_{i,j}=\bp_{12}^{-1}(\fE^\al_i)\cap \bp_{23}^{-1}(\fE^{\al+i}_j)\subset
\QL\al\times\QL{\al+i}\times\QL{\al+i+j}
$$
and
$$
\fE^\al_{j,i}=\bp_{12}^{-1}(\fE^\al_j)\cap \bp_{23}^{-1}(\fE^{\al+j}_i)\subset
\QL\al\times\QL{\al+j}\times\QL{\al+i+j},
$$
where $\bp_{ab}$ denotes the projection of the product
$\QL\al\times\QL{\al+i}\times\QL{\al+i+j}$ (resp.\
$\QL\al\times\QL{\al+j}\times\QL{\al+i+j}$) onto the product of the
$a$-th and $b$-th factors.

The definition of $\fE^\al_i$ implies that
$\fE^\al_{i,j}=\{(E_\bullet,E'_\bullet,E'''_\bullet)$ such that
$E_\bullet\supset E'_\bullet\supset E'''_\bullet$,
$E_\bullet/E'_\bullet=\CO_x[i]$ and $E'_\bullet/E'''_\bullet=\CO_y[j]$
for some $x,y\in C\,\}$.
This means that $E_\bullet/E'''_\bullet$ is an extension
of $\CO_y[j]$ by $\CO_x[i]$. But it is easy to see
that $|i-j|>1$ implies $\Ext_\RT(\CO_y[j],\CO_x[i])=0$,
hence $E_\bullet/E'''_\bullet=\CO_x[i]\oplus\CO_y[j]$.
Let $E''$ be the kernel of the composition
$E_\bullet\to\CO_x[i]\oplus\CO_y[j]\to\CO_y[j]$.
Then we have $E_\bullet\supset E''_\bullet\supset E'''_\bullet$,
$E_\bullet/E''_\bullet=\CO_y[j]$ and $E''_\bullet/E'''_\bullet=\CO_x[i]$.
This means that $(E,E',E''')\to(E,E'',E''')$ is a map from
$\fE^\al_{i,j}$ to $\fE^\al_{j,i}$. This map is certainly an isomorphism.

Now the composition of correspondences $[\fE^\al_i]\circ[\fE^\al_j]$
is the correspondence given by the cycle
$$
{\bp_{13}}_*(\bp_{12}^*[\fE^\al_i]\cap \bp_{23}^*[\fE^{\al+i}_j])=
{\bp_{13}}_*[\fE^\al_{i,j}]={\bp_{13}}_*[\fE^\al_{j,i}]=
{\bp_{13}}_*(\bp_{12}^*[\fE^\al_j]\cap \bp_{23}^*[\fE^{\al+j}_i]),
$$
i.e.\
$$
\hat e_i\hat e_j=\hat e_j\hat e_i\qquad\hat f_i\hat f_j=\hat f_j\hat f_i
$$
and Proposition follows. $\Box$

\subsection{Proposition}
\label{Serre}
Given $i,j\in I$ such that $|i-j|=1$ we have

a) $e_i^2e_j-2e_ie_je_i+e_je_i^2=0$;

b) $f_i^2f_j-2f_if_jf_i+f_jf_i^2=0$.

{\em Proof.} 
Let $j=i-1$.
Consider the space $\fE^\al_{2i+j}\subset\QL\al\times\QL{\al+2i+j}$
of all pairs
$(E,E')$ such that $E'\subset E$ and $\dim(E/E')=2i+j$. Let
$T=(0,\dots,0,T_j,T_i,0,\dots,0)=E/E'$ be the quotient representation.
Let $\bp:\fE^\al_{2i+j}\to\QL\al$ denote the map induced by the projection of
$\QL\al\times\QL{\al+2i+j}$ onto the first factor, and let
$\br:\fE^\al_{2i+j}\to C^{2i+j}$
denote the map sending a pair $(E,E')$ to 
$\sum\limits_{x\in C}\dim_x(T)x=j\,\supp T_j+i\,\supp T_i$.
Recall the diagonal stratification
$$C^{2i+j}=C^{2i+j}_{\lbr i,i,j\rbr }\sqcup C^{2i+j}_{\lbr 2i,j\rbr }\sqcup
C^{2i+j}_{\lbr i,i+j\rbr }\sqcup C^{2i+j}_{\lbr 2i+j\rbr }$$
introduced e.g. in ~\cite{k}, 1.3.

Consider the map $\bp\times \br:\fE^\al_{2i+j}\to\QL\al\times C^{2i+j}$.
This map is an isomorphism over the open set
$\CQ_\al\times C^{2i+j}_{\lbr i,i,j\rbr }\subset\QL\al\times C^{2i+j}$.
This can be proved by the same arguments as the Lemma~\ref{fibersofpr}
(the fiber over the point $(E,jx\!+\!iy\!+\!iz)$ is isomorphic to
$\PP(\Hom(E_{i-1}/E_{i-2},\CO_x))\times
\Big(\PP(\Hom(E_i/E_{i-1},\CO_y))
\times\PP(\Hom(E_i/E_{i-1},\CO_z))/\ZZZ_2\Big)$
which is a single point in our case).

On the other hand, over the subset 
$\CQ_\al\times C^{2i+j}_{\lbr i,i+j\rbr }\subset\QL\al\times C^{2i+j}$ the
fibers
of the map $\bp\times \br$ are one-dimensional (the fiber over a point
$(E,(i\!+\!j)x\!+\!iy)$ is naturally isomorphic to
$\PP(\Hom(E_i/E_{i-2},\CO_x))\times\PP(\Hom(E_i/E_{i-1},\CO_y))$ which is
$\PP^1$ in our case).
Note that for the generic element $(E,E')$ of the fiber the map
$T_{i-1}\to T_i$ in the quotient representation is non-zero.

Let $\fE^\al_{\lbr i,i,j\rbr }$ denote the closure of
$(\bp\times \br)^{-1}(\CQ_\al\times C^{2i+j}_{\lbr i,i,j\rbr })$,
and let $\fE^\al_{\lbr i,i+j\rbr }$ denote
the closure of $(\bp\times \br)^{-1}(\CQ_\al\times C^{2i+j}_{\lbr i,i+j\rbr })$.
The spaces $\fE^\al_{\lbr i,i,j\rbr }$
and $\fE^\al_{\lbr i,i+j\rbr }$ are irreducible
$(\dim\CB+2|\al|+3)$-dimensional components of $\fE^\al_{2i+j}$.

\subsubsection{Claim}\label{dim2}
All other irreducible components of $\fE^\al_{2i+j}$ have smaller dimension.

{\em Proof.}
It is just a particular case of the Proposition ~\ref{dims}. $\Box$

Now we can finish the proof of the proposition.

To this end consider the spaces
$$
{\arraycolsep=3pt
\begin{array}{lcccccccccccccc}
\fE^\al_{i,i,j} & = & \bp_{12}^{-1}(\fE^\al_i) & \cap &
\bp_{23}^{-1}(\fE^{\al+i}_i)
& \cap & \bp_{34}^{-1}(\fE^{\al+2i}_j) & \subset &
\QL\al & \times & \QL{\al+i} & \times & \QL{\al+2i} & \times & \QL{\al+2i+j}\\
\fE^\al_{i,j,i} & = & \bp_{12}^{-1}(\fE^\al_i) & \cap &
\bp_{23}^{-1}(\fE^{\al+i}_j)
& \cap & \bp_{34}^{-1}(\fE^{\al+i+j}_i) & \subset &
\QL\al & \times & \QL{\al+i} & \times & \QL{\al+i+j} & \times & \QL{\al+2i+j}\\
\fE^\al_{j,i,i} & = & \bp_{12}^{-1}(\fE^\al_j) & \cap &
\bp_{23}^{-1}(\fE^{\al+j}_i)
& \cap & \bp_{34}^{-1}(\fE^{\al+i+j}_i) & \subset &
\QL\al & \times & \QL{\al+j} & \times & \QL{\al+i+j} & \times & \QL{\al+2i+j}
\end{array}
}
$$
It is easy to see that the space $\fE^\al_{i,i,j}$
(resp.\ $\fE^\al_{i,j,i}$,
$\fE^\al_{j,i,i}$) is isomorphic to the space of triples $(E,E',F)$, where
$E\in\QL\al$, $E'\in\QL{\al+2i+j}$ such that $E'\subset E$ and 
$\dim(E/E')=2i+j$, and $F$ is a filtration (by subrepresentations) in the
quotient representation $0\subset F_1\subset F_2\subset F_3=E/E'$ of the type
$(i,i,j)$ (resp.\ $(i,j,i)$, $(j,i,i)$).

Consider the projection
$\bp_{14}:\fE^\al_{i,i,j}\to\QL\al\times\QL{\al+2i+j}$ (and two others).
It is clear that the images of $\fE^\al_{i,i,j}$ (resp.\ $\fE^\al_{i,j,i}$,
$\fE^\al_{j,i,i}$) lie in $\fE^\al_{2i+j}$ and the fibers of these projections
over the point $(E,E')$ can be identified with the set of all filtrations
$F$ in the quotient representation $T=E/E'$ of the corresponding type.

\subsubsection{Lemma}\label{filt}

a) If $\ jx\!+\!iy\!+\!iz\in C^{2i+j}_{\lbr i,i,j\rbr }$ and $T$ is
a $(2i+j)$-dimensional representation
of the quiver $\Omega$ with $\supp T_{i-1}=x$, $\supp T_i=\{y,z\}$ then
$T$ admits two filtrations of type $(i,i,j)$,
two filtrations of the type $(i,j,i)$ and two filtrations of the type $(j,i,i)$.

b) If $\ (i\!+\!j)x\!+\!iy\in C^{2i+j}_{\lbr i,i+j\rbr }$ and $T$ is 
$(2i+j)$-dimensional
representation of the quiver $\Omega$ with $\supp T_{i-1}=x$,
$\supp T_i=\{x,y\}$ and non-zero map $T_{i-1}\to T_i$ then $T$ admits
two filtrations of the type $(i,i,j)$, one filtration of the type $(i,j,i)$
and no filtrations of the type $(j,i,i)$.

{\em Proof.} Trivial. $\Box$

Now we are ready to compute the compositions of the correspondences.
\begin{multline*}
[\fE^\al_i]\circ[\fE^{\al+i}_i]\circ[\fE^{\al+2i}_j]=
(\bp_{14})_*
(\bp_{12}^*[\fE^\al_i]\cap  \bp_{23}^*[\fE^{\al+i}_i]\cap
\bp_{34}^*[\fE^{\al+2i}_j])=\\
=(\bp_{14})_*[\fE^\al_{i,i,j}]=2[\fE^\al_{\lbr i,i,j\rbr }]+
2[\fE^\al_{\lbr i,i+j\rbr }]+
\text{terms of smaller dimension.}
\end{multline*}
Similarly,
$$
[\fE^\al_i]\circ[\fE^{\al+i}_j]\circ[\fE^{\al+i+j}_i]=
2[\fE^\al_{\lbr i,i,j\rbr }]+
[\fE^\al_{\lbr i,i+j\rbr }]+\text{terms of smaller dimension,}
$$
and
$$
[\fE^\al_j]\circ[\fE^{\al+j}_i]\circ[\fE^{\al+i+j}_i]=
2[\fE^\al_{\lbr i,i,j\rbr }]+\text{terms of smaller dimension.}
$$
The Proposition in the case $j=i-1$ follows (recall Remark~\ref{rem1}).
The case $j=i+1$ can be treated similarly. $\Box$

\subsection{Proposition}
\label{h}
a) For $i\not=j$ we have $e_if_j=f_je_i$;

b) On $H^\bullet(\CQ^L_\alpha,\BQ)$ we have $e_if_i-f_ie_i=\langle
i',\alpha+2\rho\rangle$ (multiplication by a constant).

Here $i'\in X$ stands for the simple root dual to $i$, and 
$\langle\;,\,\rangle:\ X\times Y\lra\BZ$ stands for the nondegenerate pairing 
between cocharacters and weights. Finally, $2\rho\in\BN[I]$ is the sum of all 
positive coroots.

{\em Proof.} 
Consider the following spaces:
$$
{\arraycolsep=2pt
\begin{array}{lcccccccccc}
\fEF^\al_{i,j} & = & \bp_{12}^{-1}(\fE^\al_i) & \cap &
\bp_{23}^{-1}\left((\fE^{\al+i-j}_j)^T\right) & \subset &
\QL\al & \times & \QL{\al+i} & \times & \QL{\al+i-j} \\
\fFE^\al_{i,j} & = & \bp_{12}^{-1}\left((\fE^{\al-j}_j)^T\right) & \cap &
\bp_{23}^{-1}(\fE^{\al-j}_i) & \subset &
\QL\al & \times & \QL{\al-j} & \times & \QL{\al+i-j} 
\end{array}
}
$$
Here $(\fE^\al_i)^T$ denotes the subvariety in $\QL{\al+i}\times\QL\al$
transposed to $\fE^\al_i\subset\QL\al\times\QL{\al+i}$.
It is easy to see that $\fEF^\al_{i,j}$ is the space of triples
$(E,E',E''')\in\QL\al\times\QL{\al+i}\times\QL{\al+i-j}$ such that
$E\supset E'\subset E'''$ and $\fFE^\al_{i,j}$ is the space of triples
$(E,E'',E''')\in\QL\al\times\QL{\al-j}\times\QL{\al+i-j}$ such that 
$E\subset E''\supset E'''$.

Consider the projections $\bp_{13}:\fEF^\al_{i,j}\to\QL\al\times\QL{\al+i-j}$
and $\bp_{13}:\fFE^\al_{i,j}\to\QL\al\times\QL{\al+i-j}$.

Over the set $U=\{(E,E''')\ |\ E\ne E'''\}\subset\QL\al\times\QL{\al+i-j}$
(in the case $i\ne j$ we have $U=\QL\al\times\QL{\al+i-j}$)
the spaces $\fEF^\al_{i,j}$ and $\fFE^\al_{i,j}$ are isomorphic.
The isomorphisms
are given by formulas
$$
(E,E',E''')\mapsto(E,E+E''',E''')\qquad\text{and}\qquad
(E,E'',E''')\mapsto(E,E\cap E''',E''').
$$
Let $\wti \fEF^\al_{i,j}$ (resp.\ $\wti \fFE^\al_{i,j}$) denote the closure
of $\bp_{13}^{-1}(U)$ in $\fEF^\al_{i,j}$ (resp. in $\fFE^\al_{i,j}$).
In the case $i\ne j$ we have
$\wti \fEF^\al_{i,j}=\fEF^\al_{i,j}$
(resp.\ $\wti \fFE^\al_{i,j}=\fFE^\al_{i,j}$).
We have $(\bp_{13})_*[\wti \fEF^\al_{i,j}]=(\bp_{13})_*[\wti \fFE^\al_{i,j}]$.
Since
$$
[\fE^\al_i]\circ[(\fE^{\al+i-j}_j)^T]=(\bp_{13})_*[\fEF^\al_{i,j}],\qquad
[(\fE^{\al-j}_j)^T]\circ[\fE^{\al-j}_i]=(\bp_{13})_*[\fFE^\al_{i,j}],
$$
the case $i\ne j$ follows.

In the case $i=j$ it remains to compare the contribution of
components of $\fEF^\al_{i,i}$ and $\fFE^\al_{i,i}$ over the diagonal
$\QL\al @>\Delta>> \QL\al\times\QL\al$.
Let $\fEF^\al_i$ (resp.\ $\fFE^\al_i$) be the preimage
of the diagonal $\QL\al\subset\QL\al\times\QL\al$, and let
$\bp:\fEF^\al_i\to\QL\al$
(resp.\ $\barq:\fFE^\al_i\to\QL\al$) be the corresponding projection.
It is easy to see that $\fEF^\al_i$ is isomorphic
to $\fE^\al_i$ and $\fFE^\al_i$ is isomorphic to $\fE^{\al-i}_i$ (and the maps
$\bp$, $\barq$ are the same as in \ref{def}).
Hence the dimension of $\fFE^\al_i$ is equal to $\dim\CB+2|\al|-1$
which is less than the expected dimension of
$[(\fE^{\al-i}_i)^T]\circ[\fE^{\al-i}_i]$ equal to $\dim\CB+2|\al|$.
Thus the contribution of $\fFE^\al_i$ in 
$[(\fE^{\al-i}_i)^T]\circ[\fE^{\al-i}_i]$ lives in the dimension
smaller than $\dim\CB+2|\al|$.

On the other hand the dimension of $\fEF^\al_i$ is equal to $\dim\CB+|\al|+1$
which is greater than the expected dimension. According to the Intersection
Theory (see ~\cite{fu}) in this case we have
$$
[\fE^\al_i]\circ[(\fE^\al_i)^T]-[(\fE^{\al-i}_i)^T]\circ[\fE^{\al-i}_i]=
(\bp_{13})_*(c_1(\CL)),
$$
where $\CL$ is a certain line bundle on $\fEF^\al_i$ defined in ~\ref{L} below.

We know that only one component of $\fE^\al_i$ (namely $\fE^\al_{\lbr i\rbr}$)
dominates $\QL\al$. This means that
$$
[\fE^\al_i]\circ[(\fE^\al_i)^T]-[(\fE^{\al-i}_i)^T]\circ[\fE^{\al-i}_i]=
(\bp_{13})_*(c_1(\CL_{|\fE^\al_{\lbr i\rbr}}))=
\Delta_*\bp_*c_1(\CL_{|\fE^\al_{\lbr i\rbr}})+\text{terms of smaller dimension}.
$$
Since over the generic point of $\QL\al$ the fiber of $\fE^\al_{\lbr i\rbr}$
is one-dimensional we have
$$
[\fE^\al_i]\circ[(\fE^\al_i)^T]-[(\fE^{\al-i}_i)^T]\circ[\fE^{\al-i}_i]=
b^i_\al[\Delta]+\text{terms of smaller dimension,}
$$
where $b^i_\al$ is the degree of the restriction of $\CL$ to the generic
fiber of $\fE^\al_{\lbr i\rbr}$ over $\QL\al$.

Thus to prove the Proposition it suffices to compute the integers $b^i_\al$.
The calculation of $b^i_\al$ will be given in the next section. $\Box$

\subsubsection{Definition of $\CL$ and $b^i_\al$}
\label{L}
Since we are ultimately interested in the degree of $\CL$ restricted to the
generic fiber which belongs to the smooth locus of $\fE^\al_{\lbr i\rbr}$, 
below we will restrict ourselves to this smooth locus.

We have the following diagram
$$
\begin{CD}
\fE^\al_{\lbr i\rbr}      @>\id\times \bp>>   \fE^\al_i\times\QL\al     \\
@V \bp\times\id VV                           @V\id\times\id VV            \\
\QL\al\times(\fE^\al_i)^T @>\id\times\id^T>>  \QL\al\times\QL{\al+i}
\times\QL\al
\end{CD}
$$
where $\id$ denotes either identity map or natural embedding and $T$ denotes
the transposition.

According to the Intersection Theory, $\CL$ is the cokernel
of the natural map of normal bundles
$$
\CN_{\fE^\al_{\lbr i\rbr}/(\QL\al\times(\fE^\al_i)^T)}\lra
(\id\times \bp)^*\CN_{(\fE^\al_i\times\QL\al)/(\QL\al\times\QL{\al+i}
\times\QL\al)}
\lra\CL\lra0.
$$

The first term is evidently isomorphic to $\bp^*\CT_{\QL\al}$ and the second 
term is isomorphic to $\CN_{\fE^\al_{\lbr i\rbr}/(\QL\al\times\QL{\al+i})}$.
Consider the following commutative diagram
$$
\begin{CD}
                @.   \CT_{\fE^\al_{\lbr i\rbr}}  @= \CT_{\fE^\al_{\lbr
i\rbr}}\\
@.                      @VVV                    @VVV            \\
\bp^*\CT_{\QL\al} @>>>(\CT_{\QL\al\times\QL{\al+i}})_{|\fE^\al_{\lbr i\rbr}}
				@>>>                \barq^*\CT_{\QL{\al+i}} \\
@|                      @VVV                                    \\
\bp^*\CT_{\QL\al} @>>> \CN_{\fE^\al_{\lbr i\rbr}/(\QL\al\times\QL{\al+i})}
\end{CD}
$$
with exact middle row and exact middle column. This diagram implies
that we have the following exact sequence:
$$
\CT_{\fE^\al_{\lbr i\rbr}} \lra \barq^*\CT_{\QL{\al+i}} \lra \CL \lra 0.
$$
Let $D^\al_i=\barq(\fE^\al_{\lbr i\rbr})$. This is a divisor in $\QL{\al+i}$.
Let $\vphi\in\CQ_\al$, hence $\bp^{-1}(\vphi)\cong C$.
Since the restriction of $\barq$ to the open subset $\bp^{-1}(\CQ_\al)$
is an embedding we have
$$
\CL_{|\bp^{-1}(\vphi)}\cong \barq^*\CN_{D^\al_i/\QL{\al+i}}
$$
Thus we have proved the following.

\begin{lem}{balpha}
$b^i_\al=\deg \barq^*\CN_{D^\al_i/\QL{\al+i}}$

where $\vphi\in\CQ_\al$, and
$\barq:C=\bp^{-1}(\vphi)\to\QL{\al+i}$ is the map, induced
by the projection $\barq:\fE^\al_i\to\QL{\al+i}$.
\end{lem}

The calculation of these integers will be given in the next section
(see the Proposition~\ref{degN}) with the help of Kontsevich's
compactification $\CQ^K_\al$ of the space $\CQ_\al$.

\subsection{}
\label{sl}
Let us define an operator $h_i:\ H^\bullet(\CQ^L_\alpha,\BQ)\lra
H^\bullet(\CQ^L_\alpha,\BQ)$ as a scalar multiplication by
$\langle i',\alpha+2\rho\rangle$. Combining the Propositions ~\ref{triv},
~\ref{Serre}, ~\ref{h} together with the Theorem ~\ref{character} we arrive
at the following Theorem:

{\bf Theorem.} The operators $e_i,f_i,h_i(i\in I)$ extend to the action of
Lie algebra $\frak{sl}_n$ on 
$\bigoplus\limits_{\alpha\in\BN[I]}H^\bullet(\CQ^L_\alpha,\BQ)$.
The character of this $\frak{sl}_n$-module is equal to
$\dfrac{|W|e^{2\rho}}{\prod_{\theta\in R^+}(1-e^\theta)^2}$. $\Box$

\subsubsection{Remark} We would like to emphasize that the Lie algebra
$\frak{sl}_n$ acting by correspondences in the cohomology of Laumon's spaces
should be viewed as the {\em Langlands dual} of the original group
$G=SL_n$. In effect, the character of the $\frak{sl}_n$-module is naturally
a formal cocharacter of $G$ (cf. ~\ref{character}).

\section{Kontsevich's compactification}

\subsection{}
Recall the notion of a {\em stable map} (see \cite{kt} for details).

\begin{defn}{st}
A stable map is a datum $(\CC;x_1,\dots,x_m;f)$ consisting of a
connected compact reduced curve $\CC$ with $m\ge 0$ pairwise distinct marked
non-singular points and at most ordinary double singular points, and a
map $f:\CC @>>> \CX$ having no non-trivial first order infinitesimal 
automorphisms, identical on $\CX$ and $x_1,\dots,x_m$ (stability).
\end{defn}

\subsubsection{}
\begin{defn}{qk}
$\QK\al=\overline\CM_{0,0}(C\times\CB,1\oplus\al)$ is the moduli space
of stable maps to $C\times\CB$ of curves of arithmetic genus 0 with
no marked points such that 
$f_*([\CC])=1\oplus\al\in H_2(C,\ZZZ)\oplus H_2(\CB,\ZZZ)=\ZZZ\oplus\ZZZ[I]$.
\end{defn}


Given a stable map $(\CC;f)\in\QK\al$ we denote by 
$f':\CC\to C$ and $f'':\CC\to\CB$ the induced maps; we denote by $\CC_0$ the
irreducible component of $\CC$ such that $f'_*[\CC_0]=[C]$; we denote
by $\CC_1,\dots,\CC_m$ the connected components of $\CC\setminus\CC_0$;
we denote by $f_r$ (resp.\ $f'_r$, $f''_r$) the restriction of $f$
(resp.\ $f'$, $f''$) to $\CC_r$. Finally, let $\beta$ be the degree of $f''_0$ 
and $\gamma_r$ be the degree of $f_r''$ ($r=1,\dots,m$).

The space of maps $\CQ_\al$ is naturally embedded into $\QK\al$
(to every map $\vphi:C\to\CB$ we associate its graph
$\Gamma_\vphi\subset C\times\CB$) and can be identified with
the space of all stable maps $(\CC,f)$ such that $\CC$ is irreducible.
Hence we can consider $\QK\al$ as a compactification of $\CQ_\al$.

\subsection{The birational correspondence between $\QK\al$ and $\QL\al$}

Let $0\subset\CF_1\subset\dots\subset\CF_{n-1}\subset V\otimes\CO_\CB$
be the universal flag of vector bundles over the flag variety $\CB$.

A stable map $(\CC,f)\in\QK\al$ gives rise to the following flag of
vector bundles over $C$:
$$
0\subset f'_*{f''}^*\CF_1\subset\dots\subset f'_*{f''}^*\CF_{n-1}\subset
V\otimes f'_*{f''}^*\CO_\CB=V\otimes\CO_C.
$$
Note that the above inclusions are no longer inclusions of vector subbundles,
but only of coherent sheaves.
Let us denote this flag by $\Phi(\CC,f)$.

Let $U^K_\al\subset\QK\al$ denote the open subspace consisting
of all stable maps $(\CC,f)$ such that $|\ga_1+\dots+\ga_m|<2$.

\begin{lem}{Phi}
Given $(\CC,f)\in U^K_\al$ we have $\Phi(\CC,f)\in\QL\al$;
hence $\Phi$ is a map $\Phi:U^K_\al\to\QL\al$.
\end{lem}

{\em Proof.} Easy. $\Box$

\begin{rem}{remPhi}
In general, the degree of the quasiflag $\Phi(\CC,f)$ is smaller
than $\alpha$, hence the map $\Phi$ is not defined on the whole $\QK\al$.
\end{rem}

\subsubsection{}
Given a quasiflag $E_\bullet\in\QL\al$ define its graph $\Gamma_E$ in
$C\times\CB$ as follows:
\begin{multline*}
\Gamma_E=\{ (x,F_\bullet)\ |\ \text{the composition }
E_i \lra V\otimes\CO_C \lra V/F_i\otimes\CO_C\\
\text{ vanishes at the point }x\text{ for all $i=1,\dots,n-1.$}\}
\end{multline*}

Let $U^L_\al\subset\QL\al$ denote the open subspace consisting of 
all quasiflags $E_\bullet$ such that $|\deff E|<2$, where $\deff E$
stands for the defect of $E_\bullet$.

\begin{lem}{Gamma}
The graph of $E_\bullet\in U^L_\al$ is a stable curve.
Its natural embedding in $C\times\CB$ is a stable map of degree
$1\oplus\alpha$; hence the correspondence $E_\bullet\mapsto\Gamma_E$
defines a map $\Gamma: U^L_\al\to\QK\al$.
\end{lem}

{\em Proof.} Evident. $\Box$

\subsubsection{}
\begin{prop}{iso}
The maps $\Phi$ and $\Gamma$ define the mutually inverse isomorphisms
$\Phi:U^K_\al\rightleftarrows U^L_\al:\Gamma$,
which are identical on the subspace $U^K_\al\supset\CQ_\al\subset U^L_\al$.
\end{prop}

{\em Proof.} Clear. $\Box$

\subsubsection{}
Let $\tilde D^\al_i$ be closed subspace in $U^K_{\al+i}$ consisting of all
stable maps $(\CC,f)$ with $\beta=\al$, $\ga_1=i$.
\begin{lem}{D}
The maps $\Phi$ and $\Gamma$ induce the isomorphisms
$\Phi:D^\al_i\cap U^L_{\al+i}\rightleftarrows\tilde D^\al_i:\Gamma$.
Hence $\barq^*\CN_{D^\al_i/\QL{\al+i}}\cong
\tilde \barq^*\CN_{\tilde D^\al_i/\QK{\al+i}}$ where $\vphi\in\CQ_\al,\
\barq:C=\bp^{-1}(\vphi)\to\QL{\al+i}$ is the map induced by the projection
$\barq:\fE^\al_i\to\QL{\al+i}$, and $\tilde \barq$ is the composition
$C\stackrel{\barq}{\lra}\QL{\al+i}\stackrel{\Gamma}{\dasharrow}\QK{\al+i}$.
\end{lem}

{\em Proof.} Easy. $\Box$

\subsection{}\label{CBi}
Let $\CP_i\subset G$ be the minimal parabolic subgroup of type $i$ containing
the Borel subgroup $B$.  Let $\CB_i=G/\CP_i$ be the
corresponding homogenuous space, and let $\sigma_i$ stand for the natural
projection $\sigma_i:\CB=G/B\to\CB_i$. The map $\sigma_i:\CB\to\CB_i$
is a $\PP^1$-fibration, and its relative tangent bundle $\CT_{\CB/\CB_i}$
is canonically isomorphic to the line bundle $\CL_{i'}$ corresponding
to the simple root $i'$, considered as a character of $B$.

\subsubsection{Lemma}
Let $\vphi\in\CQ_\al$ be a map from $C$ to $\CB$ of degree $\al$.
The map $\tilde \barq:C\to\QK{\al+i}$ can be described as follows:
$$
x\in C\mapsto \Gamma_\vphi\cup\{x\}\times\sigma_i^{-1}(\sigma_i(\vphi(x)))
\subset C\times\CB.
$$
Here the RHS is a stable curve and $\tilde \barq(x)$ is this curve
with its natural embedding into $C\times\CB$.

{\em Proof.} Apply the definitions of $\fE^\al_i$ and of the map $\Gamma$.
$\Box$

\subsection{Proposition}\label{degN}
Given $\vphi\in\CQ_\al$ we have
$$
\deg \barq^*\CN_{D^\al_i/\QL{\al+i}}=\langle i',\al+2\rho\rangle.
$$

{\em Proof.}
Recall that the fiber of the normal bundle to the divisor $\tilde D^\al_i$
in the space of stable maps $\QK{\al+i}$ at the point 
$(\CC=\CC_0\cup\CC_1,f)$ is canonically isomorphic to
$(\CT_{\CC_0})_P\otimes(\CT_{\CC_1})_P$ where $P$ is the point
of intersection $P=\CC_0\cap\CC_1$.

{\sloppy
The canonical isomorphisms
$(\CT_{\CC_0})_P\cong{f'_0}^*(\CT_C)_{f'(P)}$,
$(\CT_{\CC_1})_P\cong{f''_1}^*(\CT_{\CB/\CB_i})_{f''(P)}$ imply
that $\tilde \barq^*\CN_{\tilde D^\al_i/\QK{\al+i}}=
\CT_C\otimes\vphi^*\CL_{i'}$.
Hence its degree equals
$$
\deg(\CT_C)+\deg\vphi^*\CL_{i'}=2+\langle i',\al\rangle=
\langle i',2\rho\rangle+\langle i',\al\rangle=\langle i',\al+2\rho\rangle.
$$
Now the Proposition follows from the Lemma~\ref{D}. $\Box$

}

\section{More on correspondences}

In this section we will follow the notations of ~\cite{lu} in the particular
case of Dynkin graph of type $A_{n-1}$.

\subsection{} We fix an orientation $\Omega=(1\lra2\lra\ldots\lra n-1)$
on the Dynkin graph with the set of vertices $I$. We fix a sequence
$\mbox{\bf{i}}\in\CX$ adapted to $\Omega$; let $\theta^1,\ldots,\theta^\nu$
be the corresponding total order on $R^+$ (see {\em loc. cit.}, \S4);
here $\nu=\frac{n(n-1)}{2}$.

G.Lusztig has introduced in {\em loc. cit.} a bijection $\bc\mapsto\Vc$
between $\BN^\nu$ and the set of isomorphism classes of representations of
the quiver $\Omega$.
For $\bc\in\BN^\nu$ we will denote by $d(\bc)\in\BN[I]$ the dimension of $\Vc$.
In the notations of ~\ref{not} we have $\bc\in\fK(\gamma)\Leftrightarrow
\gamma=d(\bc)=|\bc|$.

\subsection{}
For $\gamma\in\BN[I]$ we introduce the closed subvariety 
$\fE^\al_\gamma\subset\QL\al\times\QL{\al+\gamma}$ as follows.

{\bf Definition.} $\fE^\al_\gamma=\{(E_\bullet,E'_\bullet)$ such that
$E'_\bullet\subset E_\bullet\}$.

There are natural maps
$$
\bp:\fE^\al_\gamma\to\QL\al,\qquad \barq:\fE^\al_\gamma\to\QL{\al+\gamma}
\quad\text{and}\quad \br:\fE^\al_\gamma\to C^\gamma,
$$
where $C^\gamma$ is the configuration space (see e.g. \cite{k}, 1.3).
The first and second maps are induced by the projections of 
$\QL\al\times\QL{\al+\gamma}$ onto the first and second factors
and the third one is defined as
$$
\br((E_\bullet,E'_\bullet))=\sum_{x\in C}\dim_x(E'_\bullet/E_\bullet)\cdot x.
$$

We will be interested in irreducible components of $\fE^\al_\gamma$ of the
middle dimension $\dim\CB+2|\al|+|\gamma|$.

\subsubsection{}
Recall that for $\gamma\in\BN[I]$ we denote by $\Gamma(\gamma)$ the set
of all partitions of $\gamma$, i.e. multisubsets (subsets with multiplicities)
$\Gamma=\lbr \gamma_1,\ldots,\gamma_m\rbr $ of $\BN[I]$ with
$\sum_{r=1}^m\gamma_r=\gamma,\ \gamma_r>0$ (see e.g. ~\cite{k}, 1.3).
The diagonal stratification $C^\gamma=\sqcup_{\Gamma\in\Gamma(\gamma)}
C^\gamma_\Gamma$ was introduced e.g. in {\em loc. cit.} Recall that for
$\Gamma=\lbr \gamma_1,\ldots,\gamma_m\rbr $ we have $\dim C^\gamma_\Gamma=m$.

Given a partition $\Gamma=\lbr \gamma_1,\dots,\gamma_m\rbr \in\Gamma(\gamma)$
consider the following closed subspace of $\fE^\al_\gamma$:
$$
\fE^\al_\Gamma=\overline{(\bp\times \br)^{-1}(\CQ_\al\times C^\gamma_\Gamma)}.
$$
\begin{lem}{EG}
If $\gamma_r\in R^+$ for any $r=1,\ldots,m$ (i.e. $\Gamma\in\fK(\gamma)$
is a Kostant partition of $\gamma$),
then $\fE^\al_\Gamma$ is irreducible of dimension
$\dim\CB+2|\al|+|\gamma|$.
\end{lem}

{\em Proof.}
Since $\CQ_\al\times C^\gamma_\Gamma$ is irreducible of dimension
$\dim\CB+2|\al|+m$ ($m$ is the number of elements in
the partition), we need to check that the fibers of the projection
$\pr:(\bp\times \br)^{-1}(\CQ_\al\times C^\gamma_\Gamma)\to
\CQ_\al\times C^\gamma_\Gamma$ are irreducible of dimension
$|\gamma|-m=\sum_{r=1}^m(|\gamma_r|-1)$.
Let $\ga_r=i_{q_r}+i_{q_r+1}+\dots+i_{p_r}$. Then the fiber
over a point $(E_\bullet,\sum\ga_rx_r)\in\CQ_\al\times C^\ga_\Gamma$
is naturally isomorphic to
$$
\prod_{r=1}^m\PP(\Hom(E_{p_r}/E_{q_r-1},\CO_{x_r})).
$$
Since $E_{p_r}/E_{q_r-1}$ is locally free of rank $p_r-q_r+1=|\gamma_r|$
the Lemma follows. $\Box$

\subsubsection{Remark}\label{generic}
Denote by $\EGo\subset(\pr)^{-1}(\CQ_\al\times C^\ga_\Gamma)$ the open
subspace of $E^\al_\Gamma$ with the fiber over the point
$(E_\bullet,\sum\ga_rx_r)\in\CQ_\al\times C^\ga_\Gamma$ equal to
$$
\prod_{r=1}^m\big(\PP(\Hom(E_{p_r}/E_{q_r-1},\CO_{x_r}))\setminus
\PP(\Hom(E_{p_r}/E_{q_r},\CO_{x_r}))\big).
$$
Then for a point $(E_\bullet,E_\bullet')\in\EGo$ we have the following
decomposition of the quotient $E_\bullet/E_\bullet'$ (in the category $\RT$)
$$
E_\bullet/E_\bullet'=\bigoplus_{r=1}^m\bM_{\gamma_r}\otimes\CO_{x_r},
$$
where $\bM_{\gamma_r}$ is the indecomposable representation of
$\Omega$ (in the category of representations in vector spaces),
corresponding to coroot $\gamma_r\in R^+$.

\subsection{Proposition}\label{dims}
Dimension of any irreducible component of $\fE^\al_\ga$
is not greater than $\dim\CB+2|\al|+|\ga|$. Any component of this dimension
coincides with $\fE^\alpha_\Gamma$ for some $\Gamma\in\fK(\gamma)$
(see~\ref{EG}).

{\em Proof.}
Consider the stratification of $\QL\al\times C^\ga$ via
the defect of $E_\bullet$ at the support of $\sum\ga_rx_r\in C^\ga$,
namely
$$
\QL\al\times C^\ga=\bigsqcup
\begin{Sb}
\Gamma\in\Gamma(\ga)\\
|\ka'_1|+\dots+|\ka'_m|=\ga'\le\al
\end{Sb}
\CZ^\Gamma_{\ka'_1,\dots,\ka'_m}.
$$
Here $\CZ^\Gamma_{\ka'_1,\dots,\ka'_m}\subset\QL\al\times C^\ga$ is 
the subspace of all pairs $(E_\bullet,\sum\ga_rx_r)$ such that
$\lbr \ga_1,\dots,\ga_m\rbr =\Gamma$ and the defect of $E_\bullet$
at the point $x_r$ is of type $\ka'_r$ ($r=1,\dots,m$).


We evidently have 
\begin{multline*}
\dim\CZ^\Gamma_{\ka'_1,\dots,\ka'_m}=
\dim\CB+2\left|(\al-\sum|\ka'_r|)\right|+\sum(||\ka'_r||-K(\ka'_r))+m=\\=
\dim\CB+2|\al|-\sum||\ka'_r||-\sum K(\ka'_r)+m=
\dim\CB+2|\al|+|\ga|+\sum(1-|\ga_r+\ga'_r|-K(\ka'_r)).
\end{multline*}

\subsubsection{}
Given $(E_\bullet,\sum\ga_rx_r)\in\CZ^\Gamma_{\ka'_1,\dots,\ka'_m}$
we define $\CF(E_\bullet,\sum\ga_rx_r)$ as 
$(\pr)^{-1}((E_\bullet,\sum\ga_rx_r))$.

\begin{lem}{product}
$\CF(E_\bullet,\sum\ga_rx_r)=\prod\limits_{r=1}^m\CF(E_\bullet,\ga_rx_r)$.
\end{lem}

{\em Proof.} Absolutely similar to the proof of Proposition 2.1.2 in \cite{k}.
$\Box$

\subsubsection{}\label{hnu}
Fix $x\in C$ and $\CE_\bullet\in\QL\al$ such that the defect of
$\CE_\bullet$ at the point $x$ is of type $\ka'\in\fK(\ga')$.
Here we will study the variety $\CF(\CE_\bullet,\ga x)$.
To this end we may (and will) replace $C$ by the formal neighbourhood of $x$.

Let $\CE_q^p$ be the normalization of $\CE_q$ in $\CE_p$ and let $\ti\CE_q$ be
the normalization of $\CE_q$ in $V\otimes\CO_C$.
Then we evidently have
$\CE_q=\CE_q^q\subset\dots\subset\CE_q^{n-1}\subset\ti\CE_q$.
Given $E_\bullet\in\CF(\CE_\bullet,\ga x)$ we define 
$$
\hnu_{pq}(E_\bullet)=
\length\left(\frac{\CE_q^p\cap E_{p+1}}{\CE_q^p\cap E_p}\right)
\quad(1\le q\le p\le n-1),
$$
$$
\tnu_{pq}(E_\bullet)=
\length\left(\frac{\ti\CE_q\cap E_{p+1}}{\ti\CE_q\cap E_p}\right)
\quad(1\le q\le p\le n-1),
$$
$$
\tka_{pq}(E_\bullet)=\tnu_{pq}-\tnu_{p.q-1}\quad(1\le q\le p\le n-1)
$$
(cf. ~\cite{k} (10), (8)).
Note that $\tka$ is nothing else then the type of the defect
of $E_\bullet$, hence $\tka\in\fK(\ga'+\ga)$.

\subsubsection{Lemma}\label{less}
For all $1\le q\le p\le n-1$ we have $\hnu_{pq}\le\tnu_{pq}$.

{\em Proof.}
Since $\CE_q^p\cap E_p=\left(\CE_q^p\cap E_{p+1}\right)\bigcap
\left(\ti\CE_q\cap E_p\right)$ the Lemma follows from 
{\em loc.\ cit.}, 2.2.1. $\Box$

\subsubsection{}
Let $\fS_\hnu\subset\CF(\CE_\bullet,\ga x)$ be the subspace
of all $E_\bullet$ such that $\hnu(E_\bullet)=\hnu$.

\begin{lem}{pseu}
$\fS_\hnu$ is a pseudoaffine space of dimension
$\sum\limits_{1\le q<p\le n-1}\hnu_{pq}$.
\end{lem}

{\em Proof.} Same as the proof of {\em loc.\ cit.}, Theorem 2.3.3. $\Box$

\subsubsection{}

Since
$$
\sum_{1\le q<p\le n-1}\hnu_{pq}\le\sum_{1\le q<p\le n-1}\tnu_{pq}=
||\tka||-K(\tka)
$$
and recalling \ref{hnu}, \ref{pseu} we get the following
estimate:
$$
\dim\CF(\CE_\bullet,\ga x)\le\max_{\tka\in\fK(\ga'+\ga)}(||\tka||-K(\tka))=
|\ga+\ga'|-\min_{\tka\in\fK(\ga'+\ga)}K(\tka).
$$
Comparing it with the formula for dimension of
$\CZ^\Gamma_{\ka'_1,\dots,\ka'_m}$ we get
$$
\dim(\pr)^{-1}(\CZ^\Gamma_{\ka'_1,\dots,\ka'_m})\le
\dim\CB+2|\al|+|\ga|+
\sum_{r=1}^m\left(1-K(\ka'_r)-
\min_{\tka_r\in\fK(\ga'_r+\ga_r)}K(\tka_r)\right).
$$
Since $\ga_r\ne0$ we have $K(\tka_r)\ge1$, therefore the last
term is allways non-positive and the first part of the Proposition
follows. Furthermore, the last term is equal to zero only
if for all $r$ we have $K(\ka'_r)=0$ (hence $\ga'_r=0$) and
$\ga_r\in R^+$ for any $r$. But this is exactly
the case of Lemma~\ref{EG}. $\Box$

\subsection{}
Recall that the set $\fK(\gamma)\subset\Gamma(\gamma)$ of Kostant partitions
consists of all partitions $\lbr \ga_1,\dots,\ga_m\rbr $ of $\ga$ such that
$\gamma_r\in R^+$ for any $r=1,\ldots,m$.

We have an obvious bijection between $\fK(\gamma)$
and the set
of all $\bc\in\BN^\nu$ with $c_1\theta_1+\dots+c_\nu\theta_\nu=\ga$.

For $\bc\in\BN^\nu$ we introduce the closed subvariety $\fE_\bc^\alpha\subset
\CQ^L_\alpha\times\CQ^L_{\alpha+d(\bc)}$ as follows.

{\bf Definition.} $\fE_\bc^\alpha:=\fE^\al_\Gamma$ (see ~\ref{EG}),
where $\Gamma$ is the
partition corresponding to $\bc$.

\subsection{}
\label{pbw}
Consider the Poincar\'e dual of the fundamental class in the middle cohomology
$[\fE_\bc^\alpha]\in H^\bullet(\CQ^L_\alpha\times\CQ^L_{\alpha+d(\bc)},\BQ)$.
Viewed as a correspondence, it defines two operators:
$$e_\bc:\ H^\bullet(\CQ^L_\alpha,\BQ)\rightleftharpoons
H^\bullet(\CQ^L_{\alpha+d(\bc)},\BQ)\ :f_\bc$$
adjoint to each other with respect to Poincar\'e duality.

\subsection{Proposition}
\label{divided}
(cf. ~\cite{lu} 5.4.c) For $\bc=(c_1,\ldots,c_\nu)$ we have

a) $e_\bc=e_{\theta_1}^{(c_1)}\cdots e_{\theta_\nu}^{(c_\nu)}$;

b) $f_\bc=f_{\theta_1}^{(c_1)}\cdots f_{\theta_\nu}^{(c_\nu)}$

where $f^{(c)}$ stands for the divided power $\dfrac{f^c}{c!}$.

{\em Proof.}
Let $c_1\theta_1+\dots+c_\nu\theta_\nu=\ga$.
Let 
$$
\eac=\bp_{12}^{-1}(\fE^\al_{\theta_1})\cap\dots\cap
\bp_{N-1,N}^{-1}(\fE^{\al+\ga-\theta_\nu}_{\theta_\nu})\subset
\QL\al\times\QL{\al+\theta_1}\times\dots
\times\QL{\al+\ga-\theta_\nu}\times\QL{\al+\ga},
$$
where $N=c_1+\dots+c_\nu+1$ and $\bp_{ab}$ stands for the projection
onto the product of $a$-th and $b$-th factors. Obviously
$$
\eac=\{(E_\bullet\supset E_\bullet'\supset\dots\supset
E_\bullet^{(N)})\}\subset
\QL\al\times\QL{\al+\theta_1}\times\dots
\times\QL{\al+\ga-\theta_\nu}\times\QL{\al+\ga},
$$
hence $\bp_{1N}(\eac)\subset \fE^\al_\ga$.

This implies that
\begin{multline*}
\underbrace{[\fE^\al_{\theta_1}]\circ\dots\circ
[\fE^{\al+(c_1-1)\theta_1}_{\theta_1}]}_{c_1}
\circ\dots\circ\underbrace{[\fE^{\al+\ga-c_\nu\theta_nu}_{\theta_\nu}]
\circ\dots\circ[\fE^{\al+\ga-\theta_\nu}_{\theta_\nu}]}_{c_\nu}=
{\bp_{1N}}_*\left[\eac\right]=\\=
\sum_\Gamma a_\Gamma[\fE^\al_\Gamma]+\text{terms of smaller dimension},
\end{multline*}
where $\Gamma=\lbr \ga_1,\dots,\ga_m\rbr \in\fK(\ga)$ is a Kostant
partition of $\ga$ and $a_\Gamma$ is the
number of points in the generic fiber of $\eac$
over $\fE^\al_\Gamma$.

The space $\eac$ is naturally isomorphic to the space of triples
$(E_\bullet,E_\bullet^{(N)},F)$, where $E_\bullet\in\QL\al$,
$E_\bullet^{(N)}\in\QL{\al+\ga}$ such
that $E_\bullet^{(N)}\subset E_\bullet$ and $F$ is a filtration in the
quotient representation $0=F_N\subset\dots\subset F_1=
E_\bullet/E_\bullet^{(N)}$
of the type $(\underbrace{\theta_1,\dots,\theta_1}_{c_1},\dots,
\underbrace{\theta_\nu,\dots,\theta_\nu}_{c_\nu})$
(``of type $\bc$'' for short). Hence $a_\Gamma$ is the number
of filtrations of the type $\bc$ in the quotient $E_\bullet/E_\bullet'$ for
generic $(E_\bullet,E_\bullet')\in \fE^\al_\Gamma$.
Let $(E_\bullet,E_\bullet')\in\EGo$
(see Remark~\ref{generic}).
Then $E_\bullet/E_\bullet'=\oplus\bM_{\ga_r}\otimes\CO_{x_r}$.

Assume that $F$ is a filtration of type $\bc$ on $E_\bullet/E_\bullet'$
and let $\ti F_k=F_{c_1+\dots+c_k+1}$. Then 
$0=\ti F_\nu\subset\dots\subset\ti F_1\subset\ti F_0=E_\bullet/E_\bullet'$
(resp. $0=H^0(\ti F_\nu)\subset\dots\subset H^0(\ti F_1)\subset H^0(\ti F_0)=
H^0(E_\bullet/E_\bullet')=\bM:=\oplus\bM_{\ga_r}$) is a filtration of type
$(c_1\theta_1,\dots,c_\nu\theta_\nu)$ in the category $\RT$ (resp. in the
category of representations of $\Omega$ in vector spaces). But existence of
the latter filtration means that the isomorphism class of $\bM$ is
$\bc$ (see \cite{lu}), hence integer $a_\Gamma$ is not zero only for $\Gamma$,
corresponding to $\bc$. Therefore, in order to prove
the Proposition it remains to check that for $\Gamma$
corresponding to $\bc$ we have $a_\Gamma=c_1!\cdot\dots\cdot c_\nu!$.

By the Proposition 4.9 of {\em loc.\ cit.} the filtration $H^0(\ti F)$ in $\bM$ 
is unique, hence the filtration $\ti F$ on $E_\bullet/E_\bullet'$ is unique,
hence we need to compute the number of refinements of the filtration
$\ti F$ to a filtration $F$ of type $\bc$. But $\ti F_{k-1}/\ti F_k=
\bigoplus\limits_{t=1}^{c_k}\bM_{\theta_k}\otimes\CO_{x_{r^k_t}}$, where
$\{r^k_1,\dots,r^k_{c_k}\}=\{r\in\{1,\dots,m\}\ |\ \ga_r=\theta_k\}$. 
Hence the set of these refinements is isomorphic to the set of all orderings
of subsets $\{r^k_1,\dots,r^k_{c_k}\}$ and the Proposition follows. $\Box$

\section{Conjectures}

\subsection{}
The formal character appearing in the Theorem ~\ref{sl} is not new in the
representation theory of $\frak{sl}_n$. Let us recall its previous appearences.

First of all, let $\fn\subset\frak{sl}_n$ be the nilpotent subalgebra generated
by the simple generators $e_1,\ldots,e_{n-1}$. Let $\CN\subset\frak{sl}_n$
be the nilpotent cone. The Lie algebra $\frak{sl}_n$ acts on the cohomology
$H^\nu_\fn(\CN,\CO)$ of $\CN$ with supports in $\fn$.
The character of this module is exactly
$\dfrac{|W|e^{2\rho}}{\prod_{\theta\in R^+}(1-e^\theta)^2}$
(see e.g. ~\cite{ar}, Appendix A).

\subsection{}
Let $\zeta\in\BC$ be a root of unity of degree $p>2n$ and let $\fu$ be a
small quantum group defined by G.Lusztig for the root datum $(X,Y,\ldots)$
of type $G$ and $\zeta$ (see e.g. ~\cite{luu}).

S.Arkhipov has introduced in ~\cite{ar} the graded vector space of semiinfinite
cohomology
$H^{\frac{\infty}{2}+\bullet}_\fu$ along with the action of $\frak{sl}_n$ on it.
The (graded) character of this module is given by
$$P_G(t):=\frac{e^{2\rho}t^{-\frac{1}{2}\dim\CB}\sum_{w\in W}t^{\ell(w)}}
{\prod_{\theta\in R^+}(1-te^\theta)(1-t^{-1}e^\theta)}$$
(~\cite{ar}, Theorem 4.5).

\subsection{}
B.Feigin has conjectured (1993, unpublished) that the $\frak{sl}_n$-modules
$H^{\frac{\infty}{2}+\bullet}_\fu$ and $H^\nu_\fn(\CN,\CO)$ are isomorphic.
S.Arkhipov has checked this conjecture at the level of characters in ~\cite{ar}.

\subsection{}
\label{conjecture}
We propose the following conjecture.

{\bf Conjecture.} $\frak{sl}_n$-module
$\oplus_{\alpha\in\BN[I]}H^\bullet(\CQ^L_\alpha,\BQ)$ is isomorphic to
$H^\nu_\fn(\CN,\CO)$.

\subsubsection{Remark}
\label{freeness}
Note that
$\oplus_{\alpha\in\BN[I]}H^\bullet(\CQ^L_\alpha,\BQ)$ is evidently selfdual
(by Poincar\'e duality) while $H^\nu_\fn(\CN,\CO)$ can be easily seen to
be $\fn$-free, i.e. to posess a Verma filtration.

Thus the Conjecture admits a funny corollary that both of modules in question
are {\em tilting} (see ~\cite{ap}, chapter 1).
The conjecture would follow in turn from
this funny corollary since a tilting module is defined up to isomorphism by
its character.

\subsubsection{Remark} (V.Ostrik) Here is a sketch of a nondegenerate
$\frak{sl}_n$-invariant contragredient
self-pairing on $H^\nu_\fn(\CN,\CO)$. We will prove that
$H^\nu_\fn(\CN,\CO)$ is self-dual with respect to the standard contragredient
duality in the BGG category $\CO$. We prefer to use another duality, without
Chevalley involution on $\frak{sl}_n$, exchanging highest and lowest weight
modules. To this end it suffices to construct a nondegenerate
$\frak{sl}_n$-invariant pairing between $H^\nu_\fn(\CN,\CO)$ and
$H^\nu_{\fn_-}(\CN,\CO)$ where $\fn_-$ denotes the nilpotent subalgebra
of $\frak{sl}_n$ generated by the simple generators $f_1,\ldots,f_{n-1}$.

The desired pairing is the composition of the cup-product
$\cup:\ H^\nu_\fn(\CN,\CO)\times H^\nu_{\fn_-}(\CN,\CO)\lra
H^{2\nu}_0(\CN,\CO)$ and the trace (residue) morphism
$Res_0:\ H^{2\nu}_0(\CN,\CO)\lra H^0(0,\CO)=\BC$.
Note that the dualizing complex of $\CN$ is isomorphic to its structure
sheaf $\CO$, whence the trace morphism above.

\end{document}